\documentclass[fontsize=12pt,a4paper,headings=normal,
twoside=false,leqno,parskip=half-,abstract=true]{scrartcl}
\usepackage[english]{babel}
\usepackage[utf8]{inputenc}
\usepackage{hyperref}
\hypersetup{
 pdftitle={Hybrid Bifurcations: Periodicity from Eliminating a Line of Equilibria},
 pdfauthor={Alejandro López-Nieto, Phillipo Lappicy, Nicola Vassena, Hannes Stuke, and Jia-Yuan Dai},
 colorlinks=true,
 linkcolor=blue,
 citecolor=blue,
 filecolor=blue,
 urlcolor=blue}

\usepackage{graphicx}
\usepackage[format=plain,labelfont=bf,font=small]{caption}
\usepackage{subfigure}
\usepackage{xcolor}
\usepackage[arrow, matrix, curve]{xy}
\usepackage{tikz}
\usepackage{float}

\usepackage{tikz-cd}
\usepackage{caption}
\usepackage{amsmath,amsthm}
\usepackage{amssymb} 
\usepackage[normalem]{ulem}

\newtheorem{theorem}{Theorem}[section]
\newtheorem{lemma}[theorem]{Lemma}

\theoremstyle{definition}

\theoremstyle{remark}
\newtheorem{remark}[theorem]{\emph{Remark}}

\numberwithin{equation}{section}

\usepackage[notref,notcite,color,
final
]{showkeys}

\usepackage{cellspace}
\usepackage{diagbox}
\usepackage{pgfplots}
\usepackage{float}
\usepackage{makecell}
\usepackage{enumerate}

\newcommand\cincludegraphics[2][]{\raisebox{-0.3\height}{\includegraphics[#1]{#2}}}
\newcommand{\email}[1]{\href{mailto:#1}{#1}}

\definecolor{refkey}{rgb}{1,0,0}
\definecolor{labelkey}{rgb}{1,0,0}

\begin{document}

\title{\centering{Hybrid Bifurcations: Periodicity from Eliminating a Line of Equilibria}}

\author{
 \\
{~}\\
Alejandro López-Nieto\thanks{Department of Mathematics, National Taiwan University, No. 1, Sec. 4, Roosevelt Road, 10617 Taipei, Taiwan; \email{alopez@ntu.edu.tw}}, Phillipo Lappicy\thanks{Departamento de Análisis Matemático y Matemática Aplicada, Universidad Complutense de Madrid, Pl. de las Ciencias 3, 28040, Madrid, Spain; \email{philemos@ucm.es}}, Nicola Vassena\thanks{Interdisziplinäres Zentrum für Bioinformatik, Universität Leipzig, Härtelstraße 16-18, 04107, Leipzig, Germany; \email{nicola.vassena@uni-leipzig.de}}, 
\\
Hannes Stuke\thanks{Institut für Mathematik, Freie Universität Berlin, Arnimallee 3 D, 14195, Berlin, Germany; \email{hannes.stuke@fu-berlin.de}}, Jia-Yuan Dai\thanks{Department of Mathematics, National Tsing Hua University, No. 101, Sec. 2, Kuang-Fu Road
300 Hsinchu, Taiwan; \email{jydai@math.nthu.edu.tw}}\\
\vspace{2cm}}

\date{ }
\maketitle
\thispagestyle{empty}

\newpage

\abstract{We describe a new mechanism that triggers periodic orbits in smooth dynamical systems. To this end, we introduce the concept of hybrid bifurcations, which consists of a bifurcation without parameters and a classical bifurcation. Our main result classifies the hybrid bifurcation when a line of equilibria with an exchange point of normal stability vanishes. We showcase the efficacy of our approach by proving stable periodic coexistent solutions in an ecosystem of two competing predators with Holling's type II functional response.}


\section{Introduction}\label{Sec:1}

Classical bifurcation theory studies qualitative changes of phase portraits by tuning the value of parameters; see for example \cite{GuHo83, Vander}. In smooth dynamical systems, a bifurcation parameter $\mu \in \mathbb{R}$ serves as a stationary variable by defining $\dot{\mu}(t)=0$. This viewpoint has prompted the study of bifurcations where no variables can be made stationary, thus motivating the name: \emph{bifurcation without parameters}; see \cite{FiLi02, Li15} for details. However, the hybrid situation where a classical bifurcation occurs at a bifurcation point without parameters remains uncharted in the literature\footnote{We emphasize the different usages of the term \textit{hybrid} in the study of dynamical systems. A \textit{hybrid system} refers to a combination of continuous dynamics and discrete dynamics (see \cite{ScSc20}), while a \textit{hybrid bifurcation} consists of a bifurcation without parameters and a classical bifurcation.}. In this article, we provide the first rigorous classification of such a hybrid type of bifurcation. As a result, we prove bifurcation branches of periodic orbits and determine their local stability properties.

\begin{figure}[htbp]
    \centering
     \begin{tabular}{ c ccc c }
     \boxed{\text{Hyperbolic Hopf}} &&&&  \boxed{\text{Elliptic Hopf}}\\
      \begin{minipage}{0.30\textwidth}
      \medskip
      \cincludegraphics[width=\textwidth]{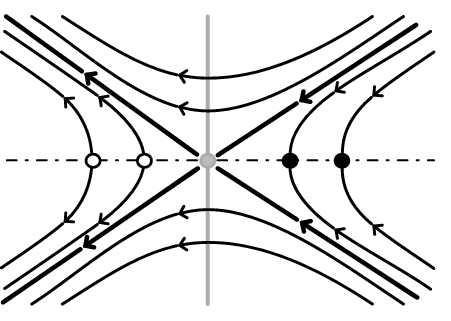}
      \end{minipage}
      &&&&
      \begin{minipage}{0.30\textwidth}
      \medskip
      \cincludegraphics[width=\textwidth]{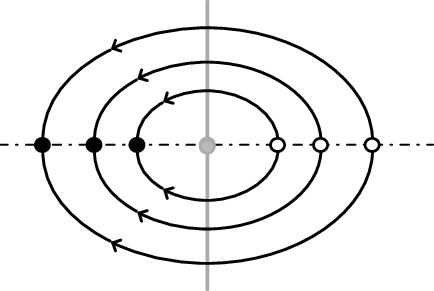}
      \end{minipage}
      \end{tabular}      
\caption{
Phase portraits in the $(r,z)$-plane of Hopf bifurcations without parameters described by the truncated cylindrical form \eqref{ToyModel}. Any foliation transverse to the line of equilibria (dash-dotted) at the Hopf point (grey circle) is also transverse to the flow nearby; thus it is not flow-invariant. The Hopf point splits the line into center-unstable equilibria (hollow circles) and center-stable equilibria (solid circles). In the hyperbolic case, two invariant manifolds (thicker lines) with rotations in $\varphi$ form an invariant cone. In the elliptic case, the Hopf point is surrounded by heteroclinic orbits. Notice that exponentially small oscillations induced by the higher-order terms not in \eqref{ToyModel} are not depicted here; see \cite[Fig. 5.2]{Li15}.
}
\label{fig1}
\end{figure}

A hybrid bifurcation consists of a bifurcation without parameters and a classical bifurcation. For simplicity of illustration, we consider ordinary differential equations (abbr., ODEs) in $(y, z; \mu) \in \mathbb{R}^n \times \mathbb{R} \times \mathbb{R}$,
\begin{equation} 
\begin{aligned} \label{first-ODEs}
\dot{y} &= f^y(y, z; \mu), \\
\dot{z} & = f^z (y,z ; \mu).
\end{aligned}
\end{equation}
Here we use the semicolon to 
distinguish between the dynamical variables $(y,z) \in \mathbb{R}^n \times \mathbb{R}$ and the stationary parameter $\mu \in \mathbb{R}$. When $\mu = 0$, we assume $f^y(0,z; 0) = 0$ and $f^z(0, z; 0) = 0$; thus $L := \{(y, z) : y = 0\}$ is a line of equilibria. Moreover, we assume that $L$ loses normal hyperbolicity at $(0,0)$. A bifurcation without parameters occurs if near $(0,0)$ there is no flow-invariant foliation transverse to $L$; see Figure \ref{fig1} for illustration. Otherwise, the $z$-equation of the ODEs \eqref{first-ODEs} is locally transformed into $\dot{z} = 0$, and thus the $z$-variable would play the role of a stationary parameter. In contrast, $\mu \neq 0$ is a stationary parameter for classical bifurcations. 

We consider the following setting for a \textit{hybrid bifurcation} in \eqref{first-ODEs}, which we illustrate in Figure \ref{fig2}:
\begin{itemize}
\item At $\mu = 0$, a line of equilibria exists and undergoes a bifurcation without parameters.
\item At $\mu \neq 0$, the line of equilibria vanishes and emanates a classical bifurcation branch.
\end{itemize}
Lines (or more generally, manifolds) of equilibria arise from symmetries of energy-based models, Hamiltonian systems, networks of oscillators, or singularly perturbed problems; see \cite{Li15}. Recent applications include electronic circuits \cite{KoSe17, Ri18} and quantum electrodynamics \cite{StGi22}. However, bifurcations without parameters require a nondegeneracy assumption (see \textbf{A4} in Section \ref{Sec:2}) and must be distinguished from the mere existence of lines of equilibria. Our approach reveals a new perspective in the following ways.
\begin{itemize}
\item We demand an exchange point of normal stability on a line of equilibria. Thus we broaden the bifurcation theory for families of solutions, which considers uniform nonhyperbolicity; see \cite{ChSb10, HaTa80}. Note that without such an exchange, bifurcations without parameters are excluded.
\item We require the elimination of the line of equilibria. In contrast, bifurcation theory without parameters only considers the scenario when the line persists; see \cite{AlFie, FiLi00, FiLi02, Li15}.
\end{itemize}

Motivated by seeking periodic orbits, we focus on Hopf bifurcations without parameters for the unperturbed case $\mu=0$; see \cite{FiLi02} and \cite[Chapter 5]{Li15}. In this situation, a line of equilibria $L$ loses normal hyperbolicity at an equilibrium $X_H$ via a pair of purely imaginary simple eigenvalues $\pm i\omega$, with $\omega > 0$. So we call $X_H$ a \textit{Hopf point}. By \cite[Theorem 5.1]{Li15}, the dynamics on a center manifold are captured by the truncated cylindrical form, 
\begin{equation}\label{ToyModel}
\begin{aligned}
    \dot\varphi&=\omega,\\
    \dot r&=rz,\\
    \dot z&=\xi r^2.
\end{aligned}
\end{equation}
Here $(\varphi, r, z) \in \mathbb{R}/2\pi\mathbb{Z} \times [0, \infty) \times \mathbb{R}$ denotes cylindrical coordinates 
\begin{equation}
(y_1,y_2,z)= X_H + (r\cos(\varphi),r\sin(\varphi),z)
\end{equation}
around the Hopf point $X_H$ and $\xi\in\{-1,1\}$ is called the \emph{discriminant}. Thus the phase portrait of \eqref{ToyModel} decomposes into a rotation $\varphi(t) = \omega t$ and planar $(r,z)$-dynamics. From \eqref{ToyModel} there are two cases depending on the discriminant $\xi \in \{-1, 1\}$:
\begin{itemize}
\item If $\xi=1$, then $X_H$ is called a \emph{hyperbolic Hopf point} and the $(r,z)$-dynamics are saddle-like with an invariant cone centered at $X_H$.
\item If $\xi=-1$, then $X_H$ is called an \emph{elliptic Hopf point}, surrounded by a continuum of heteroclinic orbits that connect equilibria in the line $L$. 
\end{itemize}
We depict the phase portraits of \eqref{ToyModel} for the two Hopf bifurcations without parameters in Figure \ref{fig1}. In contrast to classical Hopf bifurcations, no periodic orbits arise from Hopf bifurcations without parameters.

Our main conclusion is that for the perturbed case $\mu \neq 0$, under the parallel drift assumption \textbf{A5} to be defined in Section \ref{Sec:2}, the Hopf point $X_H$ undergoes a classical bifurcation to a branch of periodic orbits and we can determine their local stability properties. More precisely, up to a switching $\mu \mapsto -\mu$, we fix the direction of bifurcation for $\mu>0$ and classify the following three types of \textit{hybrid Hopf bifurcations}:
\begin{itemize}
\item \textbf{Type-H}. $X_H$ is a hyperbolic Hopf point at $\mu=0$ and emanates a branch of exponentially unstable periodic orbits for $\mu > 0$.
\item \textbf{Type-ES}. $X_H$ is an elliptic Hopf point at $\mu=0$ and emanates a branch of locally exponentially stable periodic orbits for $\mu > 0$.
\item \textbf{Type-EU}. $X_H$ is an elliptic Hopf point at $\mu=0$ and emanates a branch of exponentially unstable periodic orbits for $\mu > 0$.
\end{itemize}

Among many potential applications, we consider a system of three ODEs that describes the dynamics of two predators competing for the same prey via Holling’s type II functional response \cite{HsHuWa78, HsHuWa79, HsLiMa20}. We obtain coexistent solutions by proving a branch of stable periodic orbits from a \textbf{Type-ES} hybrid Hopf bifurcation. We emphasize that these periodic orbits do not bifurcate from boundary quadrants and thus the population size of predators can be large. Indeed, the periodic orbits obtained in the literature either require that a predator is close to extinction \cite{BuWa81, Sm82}, or rely on perturbing a conserved quantity \cite{Ke83}, or consider a singular perturbation setting \cite{LiXiYi03}. Our results overcome these limitations and provide a new mechanism that triggers stable periodic coexistence.  

This article is organized as follows: in Section \ref{Sec:2}, we state the main theorem on hybrid Hopf bifurcations and explain the core ideas of the proof. In Section \ref{Sec:3}, we review the setting of the predator-prey system with Holling’s type II functional response and state the theorem on stable periodic coexistence. At the end of the respective sections, we address the contributions of our results and further investigations. Finally, we devote Section \ref{Sec:5} to proving our results.

\section{Main result: hybrid Hopf bifurcations} \label{Sec:2}

The setting of hybrid Hopf bifurcations demands a line of equilibria $L$ that loses normal hyperbolicity at a bifurcation point $X_H$ and a stationary parameter $\mu \in \mathbb{R}$. At $\mu = 0$, the Jacobian matrix at $X_H$ has a three-dimensional center eigenspace spanned by the eigenvector along the line and two eigenvectors associated with a pair of purely imaginary simple eigenvalues. Therefore, up to center manifold reduction \cite{Vander}, the local dynamics near $X_H$ are captured by three ODEs, which we express as 
\begin{equation}\label{StandardSystem}
    \left(\begin{array}{c}
        \dot{y}\\
        \dot{z}
    \end{array}\right)=F(y,z;\mu):=
    \left(\begin{array}{c}
        f^y(y,z;\mu)\\
        f^z(y,z;\mu)
    \end{array}\right).
\end{equation}
Here $y:=(y_1,y_2)\in\mathbb{R}^2$, $z\in\mathbb{R}$, and we denote by $X := (y,z) \in \mathbb{R}^3$ the dynamical variables. 
For the sake of subsequent analysis, we assume that  $F:\mathbb{R}^2\times\mathbb{R}\times\mathbb{R}\to\mathbb{R}^3$ is sufficiently smooth. Furthermore, up to an affine linear change of coordinates, we make the following assumptions for the unperturbed case $\mu = 0$.
\begin{itemize}
\item[\textbf{A1}](\textbf{line of equilibria}). At $\mu=0$, there exists a line of equilibria $L=\{(y,z): y = 0\}$, that is, $F(0,z;0)=0$ for all $z\in\mathbb{R}$.
\item[\textbf{A2}](\textbf{spectral assumption}). The Jacobian matrix $\partial_X F(0,0;0)$ is in the Jordan form and reads
\begin{equation}\label{Jacobian}
A:= \partial_{X}F(0,0;0)=\left(
\begin{array}{ccc}
0 & -\omega &0 \\
\omega &  0 &0\\
0 & 0 & 0
\end{array}\right) \quad \mbox{for} \quad \omega > 0.
\end{equation}
\item[\textbf{A3}](\textbf{crossing assumption}). Let $\nu(z)$ denote the eigenvalues of $\partial_XF(0,z;0)$ with $\nu(0) := i\omega$. Then $\nu(z)$ crosses the imaginary axis at a nonzero speed, i.e.,
\begin{equation}
\mathrm{Re}(\nu'(0))=\partial_z(\mathrm{div}_yf^y(0,0;0))\neq 0,
\end{equation}
where $\mathrm{div}_yf^y := \partial_{y_1} f^{y_1} + \partial_{y_2}f^{y_2}$.
\item[\textbf{A4}](\textbf{nondegeneracy assumption}). $\Delta_y f^z(0,0;0) \neq 0$, where $\Delta_y f^z := \partial_{y_1}^2 f^z + \partial_{y_2}^2f^z$. 
\end{itemize}

The assumptions \textbf{A1}--\textbf{A4} ensure that the origin $(0,0) \in L$ undergoes a Hopf bifurcation without parameters; see \cite[Theorem 5.1]{Li15}. Indeed, \textbf{A1} and \textbf{A2} show that $X_H := (0,0)$ is a Hopf point. By \textbf{A3}, $X_H$ is an exchange point of normal stability on the line of equilibria $L$. The delicate assumption \textbf{A4} excludes any flow-invariant foliation near $X_H$ that is transverse to $L$, and thus the $z$-variable cannot be transformed into a stationary parameter; see the proof of \cite[Theorem 5.1]{Li15}.  

For the perturbed case $\mu \neq 0$, we eliminate the line of equilibria $L$ by introducing a parallel drift to the line at the Hopf point $X_H = (0,0)$.

\begin{itemize}
\item[\textbf{A5}](\textbf{parallel drift assumption}). $\partial_\mu f^z(0,0;0)\neq 0$ and $\partial_\mu f^y(0,0; 0) = 0$. 
\end{itemize}

Notice that \textbf{A5} excludes equilibria of the ODEs \eqref{StandardSystem} near $X_H = (0,0)$ and so eliminates the line $L$, since $f^z(0,0;\mu) \neq 0$ for sufficiently small $0 < \lvert \mu \rvert \ll 1$. Notice that $\partial_\mu f^y(0,0;0) = 0$ is merely a convenient condition for analysis, since we can always define new $y$-coordinates by the translation
\begin{equation}\label{Eq:ParamShift}    (\tilde{y}_1, \tilde{y}_2) := \left(y_1 + \frac{\mu}{\omega}\partial_\mu f^{y_2}(0,0;0),y_2-\frac{\mu}{\omega}\partial_\mu f^{y_1}(0,0;0)\right).
\end{equation}

To analyze the dynamics of the ODEs \eqref{StandardSystem} near $X_H = (0,0)$, we transform \eqref{StandardSystem} into two ODEs with a rapidly oscillating phase. Then the dynamics of \eqref{StandardSystem} are approximated by solutions of the following \textit{truncated cylindrical form}:
\begin{equation}  \label{ToyModel2}              
    \begin{aligned}
     \dot \varphi&= \omega,\\
     \dot r&=\varepsilon\beta_2rz+\varepsilon^2\left(\mu \gamma_3 r+ \mu \gamma_4 z+\beta_3r^3-(\beta_1\beta_2-\beta_4)rz^2\right),\\
     \dot z&=\varepsilon\left( \mu \gamma_5+\beta_5 r^2\right)+\varepsilon^2\left(\mu \gamma_6r+\mu(\beta_1\gamma_5+ \gamma_7)z-(\beta_1\beta_5 - \beta_6)r^2z\right).
    \end{aligned}
\end{equation}
Here the rescaling parameter $\varepsilon > 0$ extracts the leading-order terms of the expansion in the cylindrical form; see \cite[Section 7.3]{GuHo83}. In this article, we show that \eqref{ToyModel2} determines the direction of bifurcation and local stability properties of periodic orbits emerging from $X_H$. 

We highlight that the rigor of the approximation is ensured by introducing a new time $\tau$ satisfying $\dot{\tau}(t) = \dot{\varphi}(t)/\omega$ and then using averaging theory; see \cite[Chapter 6]{SanVer}. Note that the form  \eqref{ToyModel2} consists of a rapidly oscillating phase $\varphi$ and small drifting dynamics in the $(r,z)$-plane as $\varepsilon \searrow 0$. The coefficients $\beta_j$ and $\gamma_j$ will be derived in Section \ref{Sec:5}. Here, we indicate that $\beta_j$ arise from a normal form algorithm for the unperturbed case $\mu=0$, whereas $\gamma_j$ are multiplied by $\mu$ because they appear upon a $\mu$-perturbation. The explicit formulas of coefficients in terms of derivatives of $F = (f^y, f^z)$ are listed in Appendix \ref{Appendix}.

\begin{theorem}[Hybrid Hopf bifurcations]\label{T:Bifurcation}
Consider the ODEs \eqref{StandardSystem} under the assumptions \textbf{A1}--\textbf{A5}. Then there exists a $\mu_0>0$ such that \eqref{StandardSystem} possesses a local $C^1$-branch of periodic solutions $(y(t;\mu),z(t;\mu))$ for all $0<\lvert \mu \rvert <\mu_0$, satisfying
\begin{equation}\label{Eq:LimitSpace}
y(t;\mu)=O\left(\sqrt{ \lvert \mu \rvert }\right), \quad  z(t;\mu)=O(\mu) \quad \mbox{as} \quad \mu \rightarrow 0,
\end{equation}
and the direction of bifurcation is determined by
\begin{equation}\label{Eq:Direction}
\mathrm{sign}(\mu)=-\mathrm{sign}\left(\Delta_y f^z(0,0;0)\partial_\mu f^z(0,0;0)\right).
\end{equation}
Moreover, the minimal period $p(\mu)$ of the bifurcating periodic solution is continuous and satisfies
\begin{equation}\label{Eq:LimitPeriod}
\lim_{\mu\to 0}p(\mu)=\frac{2\pi}{\omega}.
\end{equation}
In addition, the discriminant defined by 
\begin{equation} \label{discriminant}
\xi:=\mathrm{sign}\left(\partial_z(\mathrm{div}_yf^y(0,0;0)) \Delta_yf^z(0,0;0)\right)
\end{equation}
satisfies $\xi \in \{-1, 1\}$ and there are three types of hybrid Hopf bifurcation:

\begin{itemize}
\item \textbf{Type-H}. If $\xi = 1$ \emph{(}i.e., $X_H = (0,0)$ is a hyperbolic Hopf point\emph{)}, then the resulting periodic orbits are exponentially unstable and possess a two-dimensional unstable manifold.
\end{itemize}
Let $\beta_j, \gamma_j$ be the coefficients of the truncated cylindrical form \eqref{ToyModel2}; see Appendix \ref{Appendix}.
\begin{itemize}
\item \textbf{Type-ES}.  If $\xi = -1$ \emph{(}i.e., $X_H$ is an elliptic Hopf point\emph{)}, then the resulting periodic orbits are locally exponentially stable if and only if
\begin{equation}\label{Eq:Stability}
2\beta_3\gamma_5^2-\beta_5\gamma_5\gamma_7+\beta_6\gamma_5^2<0.
\end{equation}
\item \textbf{Type-EU}. If $\xi=-1$, then the resulting periodic orbits are exponentially unstable with a three-dimensional unstable manifold if and only if 
\begin{equation}\label{Eq:Instability}
2\beta_3\gamma_5^2-\beta_5\gamma_5\gamma_7+\beta_6\gamma_5^2 > 0.
\end{equation}
\end{itemize}
\begin{figure}[t]
    \centering
     \begin{tabular}{c  Sc  Sc  Sc}
     &{\text{$\mu<0$}} & {\text{$\mu=0$}} & {\text{$\mu>0$}} \\ 
     \makecell{\textbf{Type-H}}&
     \begin{minipage}{0.24\textwidth}
     \cincludegraphics[width=\textwidth]{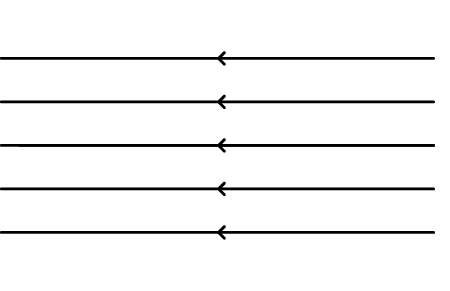}
     \end{minipage}
      & 
      \begin{minipage}{0.24\textwidth}
      \cincludegraphics[width=\textwidth]{hyphopf2.eps}
      \end{minipage}
      &
      \begin{minipage}{0.24\textwidth}
      \cincludegraphics[width=\textwidth]{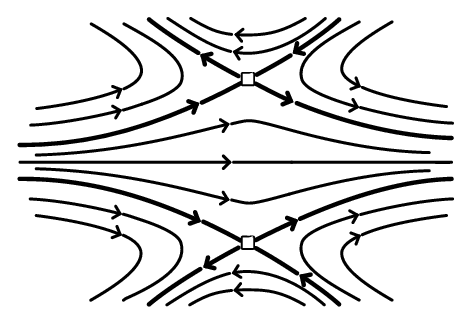}
      \end{minipage}
      \\ \hline
     \makecell{\textbf{Type-ES}}&
     \begin{minipage}{0.24\textwidth}
     \cincludegraphics[width=\textwidth]{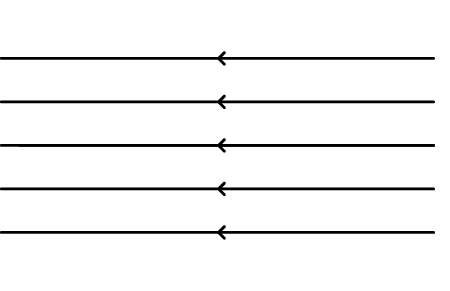}
     \end{minipage}
      &
      \begin{minipage}{0.24\textwidth}
      \cincludegraphics[width=\textwidth]{ellhopf2.eps}
      \end{minipage}
      &
      \begin{minipage}{0.24\textwidth}
      \cincludegraphics[width=\textwidth]{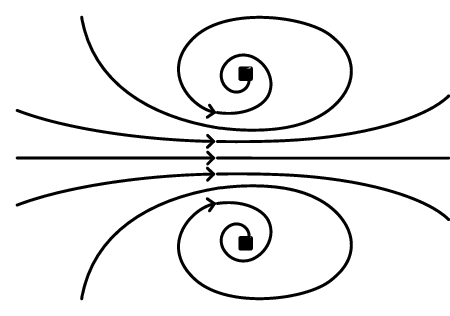}
      \end{minipage}
      \\ \hline
      \makecell{\textbf{Type-EU}}&
     \begin{minipage}{0.24\textwidth}
     \cincludegraphics[width=\textwidth]{ellhopf1.eps}
     \end{minipage}
      &
      \begin{minipage}{0.24\textwidth}
      \cincludegraphics[width=\textwidth]{ellhopf2.eps}
      \end{minipage}
      &
      \begin{minipage}{0.24\textwidth}
      \cincludegraphics[width=\textwidth]{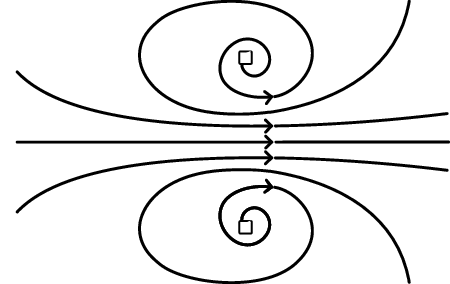}
      \end{minipage}
      \\ 
      \end{tabular}
\caption{Three types of hybrid Hopf bifurcation. At $\mu = 0$, a Hopf bifurcation without parameters occurs. The Hopf point (grey circle) splits the line of equilibria (dash-dotted) into center-unstable equilibria (hollow circles) and center-stable equilibria (solid circles). We orient $\mu$ so that the flow is parallel for $\mu<0$ and a periodic orbit (two squares, up to a phase shift $\varphi \mapsto \varphi + \pi$) bifurcates from the Hopf point for $\mu>0$. Notice that at $\mu = 0$ this orientation leads to the appearance of opposite stability properties along the line between the hyperbolic type \textbf{H} and the elliptic type \textbf{E}. In a \textbf{Type-H} hybrid Hopf bifurcation, the periodic orbit is a saddle (invariant manifolds in thick). Near an elliptic Hopf point, we can find stable periodic orbits (\textbf{Type-ES}, solid squares) or fully unstable periodic orbits (\textbf{Type-EU}, hollow squares). Notice that for $\mu > 0$, Theorem \ref{T:Bifurcation} does not exclude the existence of other invariant subsets than the periodic orbit in the region outside the parallel drifts.}
\label{fig2}
\end{figure}
\end{theorem}

\begin{remark}[Shape of bifurcation branch] \label{remark-asymptotics} The asymptotics \eqref{Eq:LimitSpace} in Theorem \ref{T:Bifurcation} shows that the branch of periodic orbits is tangent to the $y$-plane at the Hopf point $X_H$.
\end{remark}

\begin{remark}[Comparison with classical Hopf bifurcations] The application of Theorem \ref{T:Bifurcation} is comparable to that of classical Hopf bifurcation. Indeed, the assumption \textbf{A1} is included in classical Hopf bifurcations when the reference equilibrium exists for all parameter values. The assumptions \textbf{A2}--\textbf{A3} correspond to a pair of purely imaginary simple eigenvalues and the crossing condition, respectively. The main difference between Theorem \ref{T:Bifurcation} and classical Hopf bifurcations arises from the assumptions \textbf{A4--A5}. However, if we merely aim to prove the existence of periodic orbits and the direction of bifurcation, we only need to show $\Delta_y f^z(0,0; 0)$ and $\partial_\mu f^z(0,0; 0)$ are both nonzero and discuss their signs. Determining the stability in the elliptic case requires us to verify the more delicate inequalities \eqref{Eq:Stability}--\eqref{Eq:Instability}. Nevertheless, from the classical Hopf bifurcation perspective, such inequalities are analogous to determining the first Lyapunov coefficient, which is also nontrivial.
\end{remark}

We give an intuitive explanation for the formation of periodic orbits. From the truncated cylindrical form \eqref{ToyModel2}, at $\mu = 0$, the phase portrait near the line of equilibria $L = \{(\varphi, r, z) : r = 0\}$ exhibits a fast rotation $\dot{\varphi} = O(1)$ around the line $L$ and a slow drift $\dot{z} = O(\varepsilon r^2)$. In the elliptic (resp., hyperbolic) case, the drift goes from the unstable  (resp., stable) part of $L$ to the stable (resp., unstable) part of $L$; see Figure \ref{fig2} at $\mu = 0$. The perturbation $\mu > 0$ induces another slow drift $\dot{z} = O(\varepsilon \mu)$ along the already vanished line of equilibria. The appearance of periodic orbits requires balance, in the sense that the drifts have opposite directions (compare the arrows in Figure \ref{fig2} at $\mu = 0$ and at $\mu > 0$) and their amplitudes match, i.e.,  $\varepsilon r^2 = \varepsilon \mu$, which suggests the asymptotics $r = O(\sqrt{\lvert \mu \rvert})$ in \eqref{Eq:LimitSpace}. The fast rotation around $L$ motivates using averaging theory \cite[Chapter 6]{SanVer}, which is closely related to the $S^1$-equivariant normal forms in \cite{Vander} that we will adopt. 

\textbf{Ideas of proof.} The proof of hybrid Hopf bifurcations, Theorem \ref{T:Bifurcation}, is based on a rigorous derivation and analysis of the truncated cylindrical form \eqref{ToyModel2}; see Section \ref{Sec:5}. First, we adapt the normal form algorithm in \cite{Vander} to the setting under the assumptions \textbf{A1}--\textbf{A5}. Second, we find a suitable rescaling that transforms the angle variable into a fast variable. Finally, the rescaling enables averaging theory (see \cite{SanVer}) to show that an equilibrium of \eqref{ToyModel2} admits a continuation to a family of periodic orbits in the original ODEs \eqref{StandardSystem}. 

After we derive the truncated cylindrical form \eqref{ToyModel2}, its first-order truncated ODEs in $\varepsilon > 0$ read 
  \begin{equation}\label{first-order-equations}
     \begin{aligned}
     \dot{\varphi} & = \omega, \\
     \dot{r}&=\varepsilon\beta_2rz,\\
     \dot{z} &=\varepsilon(\mu\gamma_5+\beta_5 r^2),
    \end{aligned}
    \end{equation}
and have a $2\pi/\omega$-periodic solution 
\begin{equation}\label{equiequation}
    (\varphi(t),r(t),z(t))=\left(\omega t,\sqrt{\frac{-\mu \gamma_5}{\beta_5}}, 0\right) \quad \mbox{if and only if} \quad \mu \beta_5 \gamma_5 < 0,
\end{equation}
which determines the direction of bifurcation
\begin{equation} 
\mathrm{sign}(\mu) = - \mathrm{sign}(\beta_5 \gamma_5);
\end{equation}
see \eqref{Eq:Direction} and Appendix \ref{Appendix}. The linearization of \eqref{first-order-equations} at the periodic orbit in \eqref{equiequation} reads 
\begin{equation}\label{EquiLin}
    \varepsilon\left(
    \begin{array}{ccc}
        0 & 0 & 0 \\
        0 & 0 & \beta_2r_0\\
        0 & 2\beta_5r_0 & 0
    \end{array}
    \right), \quad \mbox{where   } r_0 := \sqrt{\frac{-\mu \gamma_5}{\beta_5}}.
\end{equation}
and let $A_1$ be its submatrix with respect to the $(r,z)$-variables. Then $A_1$ is invertible, which implies that the original ODEs \eqref{StandardSystem} possess a periodic solution for sufficiently small $0 < \varepsilon \ll 1$, by averaging theory; see \cite[Theorem 6.3.2]{SanVer}. For the stability analysis, we distinguish between two cases:

\begin{itemize}
\item The hyperbolic case (i.e., $\xi = 1$) corresponds to a hyperbolic matrix $A_1$ with a positive eigenvalue and a negative eigenvalue. Then the bifurcating periodic solutions belong to \textbf{Type-H} by averaging theory; see \cite[Theorem 6.3.3]{SanVer}. 
\item The elliptic case (i.e., $\xi = -1$) requires a more delicate analysis because $A_1$ has a pair of purely imaginary simple eigenvalues. Since $A_1$ is not hyperbolic, standard averaging theory does not ensure the stability properties. We overcome the obstacle by using Liouville's formula, which provides a rigorous approximation of the Floquet multipliers. Then we obtain the stability criteria \eqref{Eq:Stability}--\eqref{Eq:Instability} by analyzing the cylindrical form \eqref{ToyModel2} that involves higher-order terms. 
\end{itemize}

\textbf{Contributions and discussion.} We have established a hybrid bifurcation approach, which reveals a novel mechanism that triggers periodic orbits from eliminating a line of equilibria. We emphasize that the presence of the degeneracy involving a Hopf point obstructs the establishment of periodic orbits by well-known techniques such as classical Hopf bifurcations and geometric singular perturbation theory.

Our main result, Theorem \ref{T:Bifurcation}, classifies hybrid Hopf bifurcations and includes qualitative information such as the asymptotic behaviors of solutions \eqref{Eq:LimitSpace} and \eqref{Eq:LimitPeriod}, direction of bifurcation \eqref{Eq:Direction}, and stability criteria \eqref{Eq:Stability}--\eqref{Eq:Instability}. We indicate three remarks regarding Theorem \ref{T:Bifurcation}.

\begin{enumerate}
\item Theorem \ref{T:Bifurcation} applies to any smooth dynamical system that admits a three-dimensional center-manifold reduction (see \cite{Vander}) satisfying the assumptions \textbf{A1}--\textbf{A5}. Therefore, we have discovered a new way to obtain periodic orbits in infinite-dimensional dynamical systems governed by partial differential equations or functional differential equations.
\item The existence of periodic orbits requires the occurrence of a Hopf bifurcation without parameters. Note that classical Hopf bifurcations fulfill the assumptions \textbf{A1}--\textbf{A3}, but not the nondegeneracy assumption \textbf{A4}. Without \textbf{A4}, periodic orbits triggered by a classical Hopf bifurcation at $\mu = 0$ vanish under perturbations satisfying \textbf{A5}, since the truncated cylindrical form essentially reads
\begin{equation} \label{hopf-normal-form}
    \begin{aligned}
        \dot \varphi & = \omega, \\
        \dot r&=r(z \pm r^2),\\
        \dot z&=\mu;
    \end{aligned}
\end{equation}
also compare \eqref{hopf-normal-form} with \eqref{first-order-equations}. Indeed, at $\mu \neq 0$, the line of equilibria (i.e., the $z$-axis at $\mu = 0$) vanishes, but \eqref{hopf-normal-form} has no periodic orbits. 
\item The Jacobian matrix \eqref{Jacobian} corresponds to the fold-Hopf bifurcation of an equilibrium; see \cite{Ar94, Kuz98}. However, the presence of a line of equilibria violates the required generic condition for the standard fold-Hopf normal form method; see \cite[Lemma 8.10]{Kuz98}. Furthermore, our approach draws qualitative information for the original ODEs \eqref{StandardSystem}, rather than for the truncated versions.
\end{enumerate}

Our hybrid bifurcation approach is broadly applicable to smooth dynamical systems. The reasons are threefold: First, from a practical perspective, symmetries in smooth dynamical systems often induce a manifold of equilibria, and thus our approach is a powerful tool to find periodic solutions. Among many potential applications, we showcase the efficacy of our approach by studying the predator-prey system \eqref{full-equations}, as it possesses a symmetry of indices between the predators. Second, from a scientific perspective, real-world phenomena often present both stationary and dynamical parameters, which our hybrid bifurcation setting allows. Third, from a mathematical perspective, classical bifurcations correspond to the ODEs \eqref{first-ODEs} such that $f^z$ is identically zero, making the $z$-variable a stationary parameter. Thus they form a subset of infinite codimension in the space of nonlinearities possessing manifolds of equilibria; see \cite{Li15}. In short, classical bifurcations are a degenerate case of bifurcations without parameters. Since classical bifurcation theory is highly applicable and impactful, so are bifurcation theory without parameters and our hybrid bifurcation approach. 

\textbf{Further investigations.} We indicate four directions of research based on our hybrid bifurcation approach.

\begin{enumerate}
\item We can generate oscillatory phenomena by perturbing all known solutions triggered by Hopf bifurcations without parameters; see \cite[Chapters 6--7]{Li15}. For instance, Theorem \ref{T:Bifurcation} predicts the existence of rotating waves, which bifurcate from a family of shock waves under hyperbolic balance laws obtained in \cite{FiLi00}. The task of proving the existence is to verify the parallel drift assumption \textbf{A5}, while more delicate analysis is needed to verify the stability criteria \eqref{Eq:Stability}--\eqref{Eq:Instability}.
\item We can search dynamical systems that admit the \textit{stability boundary} between \textbf{Type-ES} and \textbf{Type-EU} hybrid Hopf bifurcations, described by  
\begin{equation}\label{Eq:Torus}
    2\beta_3\gamma_5^2-\beta_5\gamma_5\gamma_7+\beta_6\gamma_5^2=0.
\end{equation}
If \eqref{Eq:Torus} holds, then the Neimark--Sacker bifurcation (see \cite[Section 3.5]{GuHo83}) is a candidate of secondary bifurcations from periodic orbits, which yields flow-invariant tori.  

\item Existing routes to chaos in the fold-Hopf bifurcation involve homoclinic orbits of saddle-focus equilibria; see \cite[Section 7.4]{GuHo83}. In contrast, the hybrid bifurcation assumes the elimination of a line of equilibria upon a $\mu$-perturbation. Hence proving chaos in the hybrid bifurcation setting would require interactions of higher-dimensional invariant sets, and would be new to our best knowledge.

\item We can analogously consider manifolds of equilibria and obtain the hybrid counterpart of bifurcations without parameters, such as the Bogdanov--Takens point on a plane of equilibria; see \cite[Chapter 10]{Li15}. The derivation of the associated normal form follows a similar procedure to the one in Section \ref{Sec:5}. However, the truncated normal form is high-dimensional, making it challenging to completely classify the orbits that emerge through a hybrid bifurcation.
\end{enumerate}

\section{Stable periodic coexistence for competing predators}\label{Sec:3}

We consider the predator-prey system in \cite{Ke83} of two predators $x_j(t)$ for $j = 1,2$ competing for the same prey $s(t)$,
\begin{equation} \label{full-equations}
\begin{aligned}
\dot{x}_j & = \delta_j \left(\frac{s-\lambda_j}{s+\alpha_j}\right) x_j \quad \mbox{for   } j = 1, 2,
\\
\dot{s} & = s(1-s) - \frac{s}{s+\alpha_1} x_1 - \frac{s}{s+\alpha_2} x_2,
\end{aligned}
\end{equation}
with positive initial conditions $x_j(0) > 0$, $s(0) > 0$. The six parameters in the ODEs \eqref{full-equations} are positive: $\delta_j$ is the growth rate for the $j$-th predator when the ecosystem is saturated with the prey, $\alpha_j$ is the half-saturation constant for the prey, and $\lambda_j$ is the break-even concentration for the $j$-th predator because $\dot{x}_j (t) > 0$ (resp., $\dot{x}_j (t) < 0$) when $s(t) > \lambda_j$ (resp., $s(t) < \lambda_j$). We emphasize that \eqref{full-equations} is a fully rescaled system with Holling’s type II functional response; see \cite{HsHuWa78, HsHuWa79, HsLiMa20} for a derivation and interpretation of the system. Note that the carrying capacity of the prey is normalized to one. We also point out that the identical system has been studied in a chemical context \cite{BuWa81, Ke83}, where Holling’s type II functional response is known as Michaelis--Menten kinetics.

We introduce three concepts. First, we say that the $j$-th predator $x_j$ \emph{persists} if $\liminf_{t \rightarrow \infty} x_j(t) > 0$, and the persistence for the prey $s$ is defined analogously. Second, a solution of the ODEs \eqref{full-equations} is \emph{positive} if its orbit is contained in the positive octant $\mathbb{R}_+^3 := \{(x_1, x_2, s): x_1 > 0, \, x_2 > 0, \, s> 0\}$. Third, a solution is \emph{coexistent} if both predators $x_1$ and $x_2$ persist. Note that by \eqref{full-equations} the persistence of both predators implies that the prey $s$ also persists. 

We collect well-known properties from \cite{HsHuWa78} regarding the dynamics of the ODEs \eqref{full-equations}. Then we characterize the admissible parameter region for the occurrence of a hybrid Hopf bifurcation as in Theorem \ref{T:Bifurcation}. First, $\mathbb{R}_+^3$ is flow-invariant and all positive solutions are bounded in forward time. The boundary quadrants $Q_1 := \{(x_1, 0, s) : x_1 > 0, \, s > 0\}$, $Q_2 := \{(0, x_2, s) : x_2 > 0, \, s > 0\}$, and $\{(x_1, x_2, 0) : x_1 > 0, \, x_2 > 0\}$ are also flow-invariant. Second, a necessary condition for coexistent solutions is 
\begin{equation} \label{condition_lambda}
0 < \lambda_j < 1 \quad \mbox{for   } j = 1,2,
\end{equation}
i.e., the break-even concentrations of both predators are strictly less than the carrying capacity of the prey. Then there are four \textit{boundary equilibria} in $\partial \mathbb{R}_+^3$, denoted by
\begin{equation} \label{boundary-equilibria}
(0,0,0), \quad  (0,0, 1), \quad E_{1} := (x_{1}^*, 0, \lambda_1), \quad  E_{2} := (0, x_{2}^*, \lambda_2),
\end{equation}
where $x_{j}^* := (\lambda_j+\alpha_j) (1 - \lambda_j)$. 

Restricted to the invariant quadrant $Q_j$, the ODEs \eqref{full-equations} are two-dimensional and the global dynamics on $Q_j$ are well understood. More precisely, by standard linear analysis, a classical Hopf bifurcation of $E_j$ occurs on $Q_j$ when 
\begin{equation} \label{boundary-Hopf-value}
2 \lambda_j + \alpha_j - 1 = 0
\end{equation}
and triggers a \textit{boundary limit cycle} $C_j$ on $Q_j$; see \cite{Sm82}. Indeed, there is a dichotomy regarding the global dynamics: If $2\lambda_j + \alpha_j - 1 > 0$, then $E_{j}$ is globally asymptotically stable on $Q_j$; see \cite[Theorem 3.1]{ChHs98}. If $2\lambda_j + \alpha_j - 1 < 0$, then the unstable dimension of $E_{j}$ on $Q_j$ is two and the boundary limit cycle $C_j$ attracts all initial conditions in $Q_j \setminus \{E_{j}\}$; see \cite[Theorem 1]{Ch81}. 

Let $\lambda_1 := \lambda$ and $\lambda_2 := \lambda + \mu$. Our analysis begins with the unperturbed case $\mu = 0$, i.e.,  
\begin{equation} \label{equal-lambdas}
\lambda_1 = \lambda_2 = \lambda.
\end{equation}
Then both boundary equilibria $E_1$ and $E_2$ possess a center eigenvector along the following line of equilibria connecting them:
\begin{equation} \label{line}
L := \left\{(x_1, x_2, \lambda)\in \mathbb{R}_+^3 : \frac{x_1}{\lambda+\alpha_1} + \frac{x_2}{\lambda+\alpha_2} = 1 - \lambda\right\}.
\end{equation}
Notice that a line of equilibria exists if and only if $\lambda_1 = \lambda_2$, due to the ODEs \eqref{full-equations}. 

It has been shown in \cite{Wi82} that $V : \mathbb{R}^3_+ \rightarrow \mathbb{R}$ defined by
\begin{equation}\label{lyapunov}
V(x_1,x_2,s) := \frac{1}{\delta_1}\log (x_1)-\frac{\lambda + \alpha_2}{\delta_2({\lambda+ \alpha_1})}\log (x_2)
\end{equation}
has time derivative
\begin{equation} \label{time-derivative}
\dot{V}(x_1,x_2,s)=\frac{(\alpha_1-\alpha_2)(s-\lambda)^2}{(\lambda+ \alpha_1)(s+\alpha_1)(s+ \alpha_2)};
\end{equation}
thus it is a strict Lyapunov function of the ODEs \eqref{full-equations} when $\alpha_1 \neq \alpha_2$. By LaSalle's invariance principle, every positive solution converges either to the boundary or to a positive equilibrium on the line $L$. In particular, positive periodic orbits cannot exist. 

A necessary condition for periodic coexistence is thus $\lambda_1 \neq \lambda_2$, i.e., the perturbed case $\mu \neq 0$. Then the line of equilibria $L$ vanishes and only the four boundary equilibria defined in \eqref{boundary-equilibria} remain. Due to the lack of positive equilibria, no classical Hopf bifurcations occur. Nevertheless, three methods of proof have been used in the literature. First, positive periodic orbits can bifurcate from the boundary limit cycle $C_j$ via a classical transcritical bifurcation; see \cite{BuWa81, Sm82}. Such periodic orbits reside near boundary quadrants and thus the population size of one predator is small. Second, more degenerate classical bifurcations occur, for instance, by perturbing the case $\lambda_1=\lambda_2$ and $\alpha_1=\alpha_2$; see \cite{Ke83}. In this case, $V$ defined in \eqref{lyapunov} is a conserved quantity due to \eqref{time-derivative} and it can result in a tube of periodic orbits connecting both boundary quadrants $Q_1$ and $Q_2$. Then for sufficiently small $0 < \lambda_2 - \lambda_1 \ll 1$ and $0 < \alpha_2 - \alpha_1 \ll 1$, the tube breaks and a classical bifurcation yields stable positive periodic orbits. Third, geometric singular perturbation theory is applicable, for instance, for sufficiently small growth rates $0 < \delta_1, \delta_2 \ll 1$; see \cite{LiXiYi03}. We emphasize that neither the conserved quantity nor the singularly perturbed setting mentioned above admits bifurcations without parameters; see \cite[Chapter 1]{Li15}.

Our goal is to obtain stable positive periodic orbits via a hybrid Hopf bifurcation as in Theorem \ref{T:Bifurcation}. Since the ODEs \eqref{full-equations} remain unchanged after a relabeling of indices, without loss of generality we fix the order
\begin{equation}\label{order}
    0<\alpha_1\leq\alpha_2<1,
\end{equation}
where the constraint $\alpha_2<1$ follows from the ecological assumptions used in the unrescaled system; see \cite{Ke83}. Then the line of equilibria $L$ defined in \eqref{line} loses normal hyperbolicity only if we consider
\begin{equation} \label{order-2}
1-\alpha_2<2\lambda<1-\alpha_1.
\end{equation}
Under the constraint \eqref{order-2}, $L$ loses normal hyperbolicity only at a Hopf point, given by 
\begin{equation}\label{HopfPoint}
X_H := \left(\frac{(\lambda + \alpha_1)^2(2\lambda + \alpha_2 - 1)}{\alpha_2 - \alpha_1}, \frac{(\lambda + \alpha_2)^2(1 - 2\lambda- \alpha_1)}{\alpha_2 - \alpha_1}, \lambda\right).
\end{equation}

Combining the constraints \eqref{condition_lambda}, \eqref{order}, and \eqref{order-2}, we seek periodic orbits bifurcating from the Hopf point $X_H$ within the \textit{admissible parameter region} determined by 
\begin{equation}\label{ParRegion}
\delta_1, \delta_2>0,\quad 0<\lambda<\frac{1}{2},\quad 0<\alpha_1<1-2\lambda,\quad \text{and}\quad 1-2\lambda<\alpha_2<1.
\end{equation}

\begin{remark} \label{remark-tetrahedra}
We view the parameter region \eqref{ParRegion} as a $(\delta_1, \delta_2)$-family of tetrahedra in the $(\lambda, \alpha_1, \alpha_2)$-space, for the sake of the stability analysis in the proof of Theorem \ref{T:StableCoexistence}.
\end{remark}

We first show that only elliptic Hopf bifurcations occur in the admissible parameter region \eqref{ParRegion}. Hence the predator-prey system \eqref{full-equations} excludes \textbf{Type-H} hybrid Hopf bifurcations.

\begin{lemma}[All Hopf points are elliptic] \label{lemma-all-elliptic}
Consider the unperturbed case $\lambda_1 = \lambda_2 =: \lambda$ and the ODEs \eqref{full-equations} within the admissible parameter region \eqref{ParRegion}. Then $X_H$ is an elliptic Hopf point.
\end{lemma}

Lemma \ref{lemma-all-elliptic} ensures that the Hopf point $X_H$ is locally surrounded by heteroclinic orbits that connect equilibria on the line $L$; see \cite[Theorem 5.1]{Li15}. By Theorem \ref{T:Bifurcation}, we expect that eliminating the line of equilibria $L$ leads to a hybrid Hopf bifurcation of elliptic type. Remarkably, only \textbf{Type-ES} occurs. 

\begin{theorem}[Stable periodic coexistence]\label{T:StableCoexistence}
Consider the perturbed case $\lambda_1 = \lambda$, $\lambda_2 = \lambda + \mu$, and a fixed choice of parameters in the admissible region \eqref{ParRegion}. Then the ODEs
     \begin{equation} \label{unfolding}
        \begin{aligned}
             \dot{x}_1 & = \delta_1 \left(\frac{s-\lambda}{s+\alpha_1}\right) x_1,\\
             \dot{x}_2 & = \delta_2 \left(\frac{s-(\lambda+\mu)}{s+\alpha_2}\right) x_2,\\
            \dot{s} & = s(1-s) - \frac{s}{s+\alpha_1} x_1 - \frac{s}{s+\alpha_2} x_2,
        \end{aligned}
    \end{equation}
undergo a \textbf{Type-ES} hybrid Hopf bifurcation from the Hopf point $X_H$. Moreover, the bifurcation branch appears for $\mu>0$. Consequently, there exists a $\mu_0 > 0$ such that the ODEs \eqref{unfolding} possess a stable positive periodic solution for all $0 < \mu < \mu_0$.
\end{theorem}

\textbf{Ideas of proof.} We prove Lemma \ref{lemma-all-elliptic} and Theorem \ref{T:StableCoexistence} by computing the coefficients $\beta_j$ and $\gamma_j$ in Theorem \ref{T:Bifurcation} and verifying the stability inequality \eqref{Eq:Stability}; see Section \ref{Sec:5}.

\textbf{Contributions and discussion.} We have successfully applied Theorem \ref{T:Bifurcation} to obtain stable positive periodic orbits of the predator-prey system \eqref{full-equations}. Compared to the relevant literature \cite{BuWa81, Ke83, LiXiYi03}, our hybrid bifurcation approach relaxes the constraint on parameters, in the sense that $0 < \alpha_2 - \alpha_1 \ll 1$ or $0 < \delta_1, \delta_2 \ll 1$ are not assumed. Moreover, the bifurcating periodic orbits can reside far from the boundary quadrants, because the Hopf point can be located at any point on the line of equilibria.

We stress two advantages of Theorem \ref{T:StableCoexistence} from an ecological viewpoint. First, the admissible parameter region \eqref{ParRegion} is an observation on the simpler boundary dynamics. Indeed, the inequalities $1-\alpha_2<2\lambda<1-\alpha_1$ are equivalent to the instability of the boundary equilibrium $E_1$ and the stability of another boundary equilibrium $E_2$; see \eqref{boundary-equilibria}. Second, to obtain stable periodic coexistence, we only need to estimate the ratio $\lambda_1/\lambda_2$ between the break-even concentrations, rather than measuring the values of all parameters involved in the system \eqref{full-equations}.

\begin{figure}[htbp]
    \centering
      \cincludegraphics[width=0.45\textwidth]{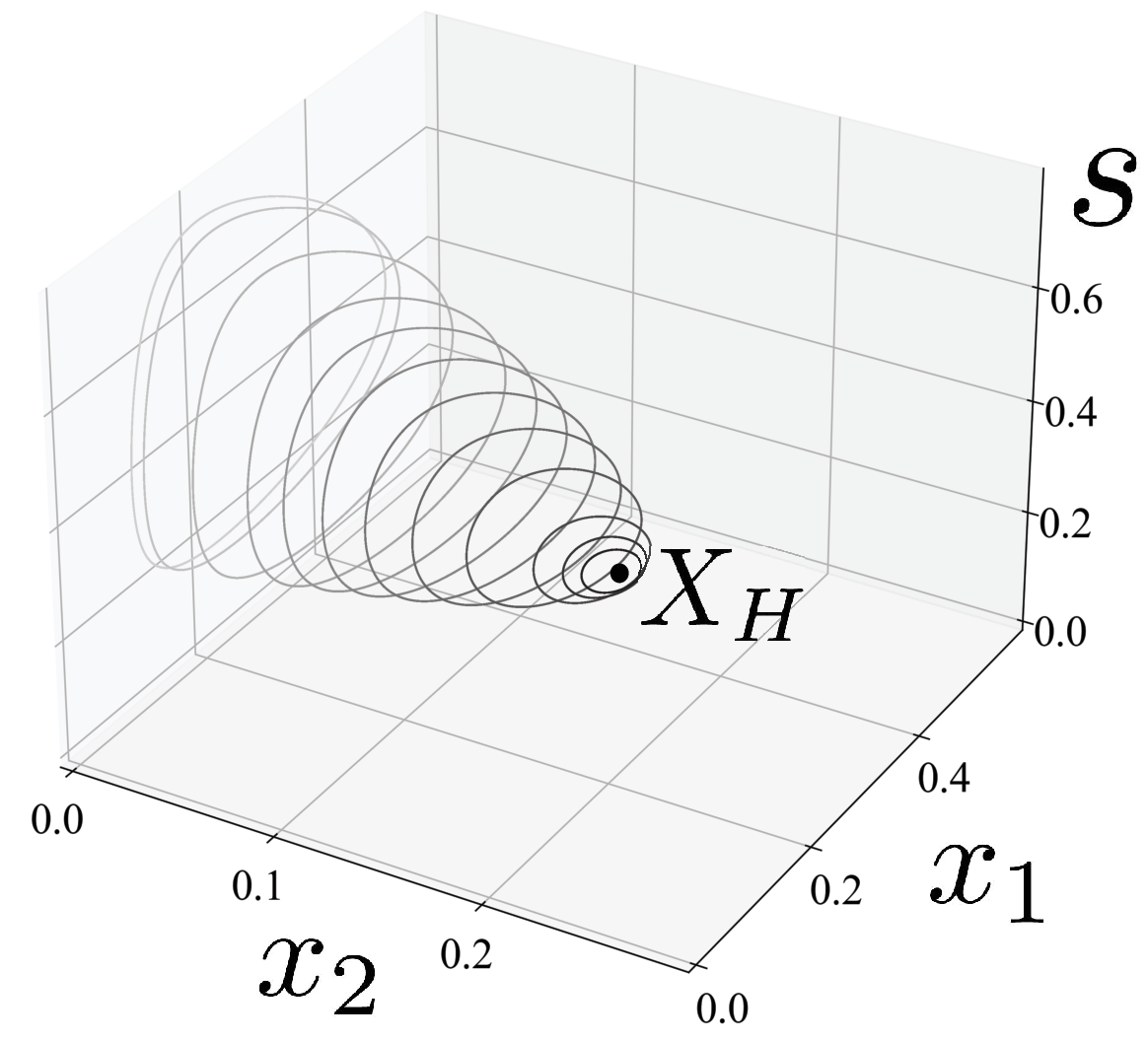}
\caption{
Parameter continuation of periodic orbits that emerge from the Hopf point towards the boundary plane $Q_1$. The parameter values are $\delta_1=0.8$, $\delta_2=0.5$, $\alpha_1=0.1$, $\alpha_2=0.2$, $\lambda=0.4$. The Hopf point, with $\mu=0$, is given by $X_H=(0.2133,0.1667,0.4)$. The twelve periodic orbits are obtained from smaller (darker) to larger (lighter) by increasing $\mu > 0$. Our choice for $\mu$ to depict such global continuation is, respectively,  $\mu =  0.0005, 0.001, 0.002,  0.005, 0.008, 0.011, 0.014, 0.017, 0.02, 0.025, 0.03, 0.05$.}

\label{figCONT}
\end{figure}

\textbf{Further investigations.} 
Supported by numerical evidence, we conjecture that the bifurcation branch of periodic orbits, obtained in Theorem \ref{T:StableCoexistence}, connects through a $\mu$-parameter continuation to the boundary limit cycle $C_1$ on $Q_1 = \{(x_1, 0, s): x_1 > 0, \, s > 0\}$ for $\mu>0$; see Figure \ref{figCONT}.

\section{Proofs of the main results}\label{Sec:5}

In this section, we prove the main results: Theorem \ref{T:Bifurcation} on hybrid Hopf bifurcations and Theorem \ref{T:StableCoexistence} on stable periodic coexistence for competing predators. Our proof of Theorem \ref{T:Bifurcation} requires a derivation of the cylindrical forms near the Hopf point at $\mu = 0$ (see Lemma \ref{L:UnpertNormalForm}) and also at $\mu \neq 0$ (see Lemma \ref{L:PertNormalform}). Then we prove Theorem \ref{T:Bifurcation} by analyzing the truncated cylindrical forms in suitably rescaled variables. Finally, we prove Theorem \ref{T:StableCoexistence} by computing the coefficients in the cylindrical form for the predator-prey system \eqref{full-equations} and verifying the stability inequality \eqref{Eq:Stability}.

We begin with $\mu = 0$ and introduce the following cylinder punctured by the line $L = \{(y,z) : y = 0\}$:
\begin{equation}
    \Omega_\rho:=\left\{(y,z)\in\mathbb{R}^2\times\mathbb{R} : 0<\lvert y \rvert<\rho\text{ and } \lvert z\rvert <\rho\right\}.
\end{equation}

\begin{lemma}[Unperturbed cylindrical form]\label{L:UnpertNormalForm}
Consider the ODEs \eqref{StandardSystem} with $\mu = 0$ and the assumptions \textbf{A1}--\textbf{A4}. Then there exist a $\rho>0$ and a diffeomorphism 
\begin{equation}
\Phi:\Omega_\rho\to \mathbb{R}/2\pi\mathbb{Z}\times (0,\rho)\times (-\rho,\rho), \quad \Phi(y,z) =: (\varphi, r, z)
\end{equation} 
such that \eqref{StandardSystem} takes the following unperturbed cylindrical form as $\lvert (r,z)\rvert \rightarrow 0$: 
\begin{equation}\label{cylcoor}
\begin{aligned}
\dot{\varphi}&=\omega+\omega\beta_1 z +O(\lvert (r,z) \rvert^2),\\
\dot{r}&=\beta_2rz+\beta_3r^3+\beta_4rz^2+O(\lvert (r,z) \rvert^4),\\
\dot{z}&=\beta_5 r^2+\beta_6r^2z+O(\lvert (r,z)\rvert^4),
\end{aligned}
\end{equation}
for $(\varphi,r,z)\in \mathbb{R}/2\pi\mathbb{Z}\times (0,\rho)\times (-\rho,\rho)$. See Appendix \ref{Appendix} for the explicit formulas of the coefficients $\beta_j$.
\end{lemma}
\begin{remark}
We highlight that the higher-order terms in the unperturbed cylindrical form \eqref{cylcoor} may depend on the angle variable $\varphi$. However, since $\varphi$ is defined on a compact domain $\mathbb{R}/2\pi\mathbb{Z}$, we omit the $\varphi$-dependence in the Landau notation.
\end{remark}
\begin{proof}[\textbf{Proof of Lemma \ref{L:UnpertNormalForm}.}]
Since $\mu = 0$ is fixed, all functions are considered to be $\mu$-independent for simplicity of notation. We follow \cite{Vander} to obtain a near-identity transformation $(y,z)=\Psi(\tilde{y},\Tilde{z})$ around the origin, which transforms the ODEs \eqref{StandardSystem} into their normal form. Omitting the tilde in the new coordinates, \eqref{StandardSystem} reads 
    \begin{equation}\label{normalform}
        \left(\begin{array}{c}
        \dot{y}\\
        \dot{z}
        \end{array}\right)
        =G(y,z):=
        \left(\begin{array}{c}
        g^y(y,z)+\eta^y(y,z)\\
        g^z(y,z)+\eta^z(y,z)
        \end{array}\right).
    \end{equation}
Here $(g^y, g^z) :=T^m G$ denotes the $m$-th Taylor polynomial of $G$ at $(y,z) = (0,0)$ and thus $(\eta^y, \eta^z) :=G-T^mG$ is the Taylor remainder of order $O(\lvert(y,z)\rvert^{m+1})$. Following \cite[Corollary 4.5]{Vander}, we define the operators $P, K : C^0 \rightarrow C^0$ by
    \begin{align}\label{averaging}
        (P F)(x) & :=\frac{\omega}{2\pi}\int_{0}^{\frac{2\pi}{\omega}}\mathrm{e}^{-A \theta}F(\mathrm{e}^{A\theta}x)\,\mathrm{d}\theta, \\
        (KF)(x)& :=\frac{\omega}{2\pi}\int_{0}^{\frac{2\pi}{\omega}}\left(\int_{0}^{\sigma}\mathrm{e}^{-A\theta}F(\mathrm{e}^{A\theta}x)\,\mathrm{d}\theta\right)\,\mathrm{d}\sigma,
    \end{align}
and the normal form algorithm proceeds recursively as follows:
\begin{itemize}
\item Set $G_1:=T^1G=\mathrm{A}$ and $\Psi_1:=T^1\Psi= I$, where $I : \mathbb{R}^3 \rightarrow \mathbb{R}^3$ denotes the identity matrix.
\item We denote $\tilde{T}^m:=T^m-T^{m-1}$ and compute the $m$-th Taylor terms $G_m:=\tilde{T}^mG$ and $\Psi_m:=\tilde{T}^m\Psi$ recursively from the $(m-1)$-th Taylor terms by applying the formula
\begin{equation}\label{algorithm}
\begin{aligned}     G_m&:=P\left(\Tilde{T}^m\left((T^mF-A)\circ T^{m-1}\Psi-D(T^{m-1}\Psi-I)\circ(T^{m-1}G-A)\right)\right),\\
\Psi_m&:=-K\left(\Tilde{T}^m\left((T^mF-A)\circ T^{m-1}\Psi-D(T^{m-1}\Psi-I)\circ(T^{m-1}G-A)\right)\right).
\end{aligned}
\end{equation}
\end{itemize}
Here the circle denotes function composition.

The normal form theorem \cite[Theorem 2.1]{Vander} guarantees that for each $m \in \mathbb{N}\cup \{0\}$ the $m$-th Taylor polynomial $T^mG$ is equivariant with respect to rotations around the $z$-axis, thus motivating cylindrical coordinates
\begin{equation} \label{cylindrical-coordinates}
(y_1,y_2,z)=(r\cos(\varphi),r\sin(\varphi),z).
\end{equation}
Notice that the punctured cylinder $\Omega_\rho$ is chosen as a domain because cylindrical coordinates \eqref{cylindrical-coordinates} are singular at $r = 0$, i.e., the $z$-axis. Then we take sufficiently small $0 < \rho \ll 1$ so that the normal form algorithm \eqref{algorithm} is applicable up to $m=3$ in the closure of $\Omega_\rho$. By transforming the ODEs \eqref{normalform} into cylindrical coordinates, we derive \eqref{cylcoor}, where the $\varphi$-equation is truncated to the second order because the higher-order terms play no role in subsequent analysis.  Then the diffeomorphism $\Phi$ is the near-identity transformation $\Psi$ in cylindrical coordinates \eqref{cylindrical-coordinates}.
\end{proof}

Next, we introduce a $\mu$-perturbation satisfying the parallel drift assumption \textbf{A5}. 

\begin{lemma}[Perturbed cylindrical form]\label{L:PertNormalform}
Consider the ODEs \eqref{StandardSystem} with $\mu \neq 0$ and the assumptions \textbf{A1}--\textbf{A5}. Let $U_\rho:=\{(y,z)\in \Omega_\rho : \lvert z \rvert< \lvert y \rvert \}$. Then the restriction of the diffeomorphism $\Phi$ to $U_\rho$ transforms \eqref{StandardSystem} into the following perturbed cylindrical form as both $r, \mu \rightarrow 0$: 
\begin{equation}\label{perturbedform}
\begin{aligned}
\dot{\varphi}&=\omega+\mu\gamma_1(\varphi) + \omega\beta_1 z +\mu\gamma_2(\varphi)\frac{z}{r}+  \mu^2 O(r^{-1}) + \mu O(r)+O(r^2),\\
\dot{r}&= \mu\gamma_3(\varphi) r+\mu \gamma_4(\varphi) z + \beta_2rz+\beta_3r^3+ \beta_4rz^2  + \mu^2 O(1) + \mu O(r^2)+O(r^4),\\
\dot{z}&=\mu\gamma_5+ \mu\gamma_6(\varphi)r+\mu \gamma_7z + \beta_5 r^2+\beta_6r^2z+ \mu^2 O(1) + \mu O(r^2)+O(r^4),
    \end{aligned}
    \end{equation}
for $(\varphi,r,z)\in \Phi(U_\rho) = \mathbb{R}/2\pi\mathbb{Z}\times \{(r,z) : \lvert z \rvert<r, 0<r<\rho\}$. See Appendix \ref{Appendix} for the formulas of the coefficients $\beta_j$ and $\gamma_j$.
\end{lemma}

\begin{remark}
Notice that $\gamma_j$ with $j = 1,2,3,4,6$ are functions of the angle variable $\varphi$, but $\gamma_5$ and $\gamma_7$ are constants. The reason is that $\gamma_5$ and $\gamma_7$ appear as a perturbation in the $z$-direction and are not involved in the transformation to cylindrical coordinates \eqref{cylindrical-coordinates} of the perturbed normal form \eqref{perturbedform}.
\end{remark}
\begin{remark} \label{remark-error-terms}
We factor out $\mu$ in the Landau notation in \eqref{perturbedform} to highlight that the $\mu$-dependent error terms are at least first and second order in $\mu$. These terms are needed for deriving the correct error bounds upon the rescaling \eqref{perturbedform}; see the proof of Theorem \ref{T:Bifurcation}.
\end{remark}

\begin{remark}[$\rho$-independence of the perturbed cylindrical form]\label{remark-radius-independence}
Notice that the formulas in Appendix \ref{Appendix} depend on the nonlinearity $f$, only. Hence the perturbed cylindrical form \eqref{perturbedform} is independent of the choice of a radius $\rho > 0$, provided that $\rho$ is sufficiently small. This independence will be crucial in our proof of Theorem \ref{T:Bifurcation}, where we will fix $\rho$ a posteriori.
\end{remark}

\begin{proof}[\textbf{Proof of Lemma \ref{L:PertNormalform}.}]
By expanding the vector field $F$ with respect to $\mu$ and using the parallel drift  assumption \textbf{A5}, we obtain 
\begin{equation}
\begin{aligned}
F(y,z;\mu) &=F(y,z;0)+\mu\left(\begin{array}{c}
             0  \\
             \gamma_5 
        \end{array}\right)+\mu \partial_\mu \partial_XF(0,0;0)\left(\begin{array}{c}
             y \\
             z
        \end{array}\right)+\mu O(\lvert(y,z)^2\rvert),
\end{aligned}
\end{equation}
where $\gamma_5:=\partial_\mu f^z(0,0;0)$. We apply the change of coordinates  $(y,z) = \Psi(\tilde{y},\tilde{z}) = \Psi_1(\tilde{y},\Tilde{z}) + \Psi_2(\tilde{y},\tilde{z}) + O(\lvert (\tilde{y}, \tilde{z})\rvert^3)$ as in the proof of Lemma \ref{L:UnpertNormalForm}. Recall that $\Psi_1=T^1\Psi = I$, i.e., $\Psi$ leaves the linear terms in $(y,z)$ unchanged. Hence $M:=I-D\Psi_2(\tilde{y}, \tilde{z})$ approximates the inverse of the Jacobian matrix, i.e.,
\begin{equation}
\begin{aligned}
    \|D\Psi(\tilde{y}, \tilde{z})^{-1} - M\| = O(\lvert(\tilde{y}, \tilde{z})\rvert^2).
\end{aligned}
\end{equation}
Next we denote 
\begin{equation}
    F_1 := \left(\begin{array}{c}
             0  \\
             \gamma_5 
        \end{array}\right) \in \mathbb{R}^3\quad\text{and}\quad F_2(y,z):= \partial_\mu \partial_XF(0,0;0)\left(\begin{array}{c}
             y \\
             z
             \end{array}\right),
\end{equation}
so that, omitting the tilde, we obtain the perturbed normal form
    \begin{equation}\label{normalform2}
    \begin{aligned}
    \left(
        \begin{array}{c}
             \dot y\\
             \dot z
        \end{array}
    \right) &= D\Psi(y,z)^{-1} F(\Psi(y,z); \mu)\\
    &= D\Psi(y,z)^{-1} F(\Psi(y,z); 0) + \mu MF_1 + \mu F_2(\Psi(y,z)) + \mu O(\lvert (y,z)^2\rvert)\\
    &= 
    \left(
    \begin{array}{c}
         g^y(y,z)  \\
         g^z(y,z)
    \end{array}\right) + \mu MF_1 + \mu F_2(\Psi(y,z)) +\mu O(\lvert(y, z)\rvert^2) + O(\lvert (y, z)\rvert^3).
    \end{aligned}
    \end{equation}
    Here $g^y$, $g^z$ are given in \eqref{normalform} and $O(\lvert(y,z)\rvert^3)$ is the $\mu$-independent error term corresponding to the functions $\eta^y$ and $\eta^z$ in \eqref{normalform}. By transforming \eqref{normalform2} into cylindrical coordinates \eqref{cylindrical-coordinates}, we obtain \eqref{perturbedform}. 
    Notice that, for $m\geq 2$, terms of the form $z^m/r$ in the $\varphi$-equation of \eqref{perturbedform} remain bounded in the domain $\Phi(U_\rho)$ and belong to $\mu O(r)$. Moreover, due to the parallel drift assumption \textbf{A5}, terms of the order $1/r$ only appear in the $\varphi$-equation and belong to $\mu^2 O(r^{-1})$.
    Finally, we compute the coefficients $\gamma_j$ using computer algebra software and keep track of the $\mu$-dependent term $\mu MF_1 + \mu F_2$ as we transform \eqref{normalform2} into cylindrical coordinates.
\end{proof}

\begin{proof}[\textbf{Proof of Theorem \ref{T:Bifurcation}}.] We split the proof into four steps. First, we rescale the perturbed cylindrical form \eqref{perturbedform} to extract the leading-order terms. Second, we introduce a new time variable that allows us to integrate the $\varphi$-equation explicitly. Third, we prove the existence of a bifurcation branch of periodic orbits by studying the first-order truncated ODEs of \eqref{perturbedform}. Fourth, we determine the local stability properties of the bifurcating periodic orbits by approximating the associated Floquet multipliers. 

\textbf{Step 1 (rescaling):} We simplify the perturbed cylindrical form \eqref{perturbedform} by introducing the rescaled variables 
    \begin{equation}\label{kappaScaling}
        r=:\varepsilon\Tilde{r},\quad z=:\varepsilon\Tilde{z},\quad \mu=:\varepsilon^2\Tilde{\mu} \quad \mbox{with} \quad \varepsilon > 0.
    \end{equation}
By Remark \ref{remark-radius-independence}, the validity of \eqref{perturbedform} is independent of a radius $\rho > 0$. Hence we choose $\rho = \varepsilon$ without loss of generality. Substituting \eqref{kappaScaling} into \eqref{perturbedform}, we obtain
    \begin{equation}\label{perturbedform2}
    \begin{aligned}
       \dot{\varphi}&=\omega + \varepsilon\omega\beta_1 \tilde{z} + O(\varepsilon^2),\\      \dot{\tilde{r}}&=\varepsilon\beta_2\tilde{r}\tilde{z} + \varepsilon^2\left(\tilde{\mu}\gamma_3(\varphi) \tilde{r} + \tilde{\mu}\gamma_4(\varphi) \tilde{z} + \beta_3 \tilde{r}^3 + \beta_4\tilde{r}\tilde{z}^2\right) +O(\varepsilon^3),\\
        \dot{\tilde{z}}&=\varepsilon\left(\tilde{\mu}\gamma_5 + \beta_5 \tilde{r}^2\right) + \varepsilon^2\left(\tilde{\mu}\gamma_6(\varphi)\tilde{r} + \tilde{\mu}\gamma_7\tilde{z} + \beta_6\tilde{r}^2\tilde{z}\right) + O(\varepsilon^3),
    \end{aligned}
    \end{equation}
where the rescaled variables satisfy $(\varphi,\tilde{r},\tilde{z})\in \Phi(U_{\rho / \varepsilon}) = \Phi(U_{1})$ and $\tilde{\mu}\in \mathbb{R}$. Notice that the term $1/r$ in the $\varphi$-equation of \eqref{perturbedform} is scaled out.

\textbf{Step 2 (decoupling in new time):} 
Since $\omega > 0$ is given, we choose a sufficiently small $0 < \varepsilon_0 \ll 1$ such that the right-hand side of the $\varphi$-equation in the ODEs \eqref{perturbedform2} is positive within $\Phi(U_{1})$ for all $0 < \varepsilon < \varepsilon_0$. Then we define a new time $\tau = \tau(t)$ by solving

\begin{equation}\label{Eq:TimeRescale}
    \dot{\tau}(t) = \frac{\dot{\varphi}(t)}{\omega}.
\end{equation}
Throughout the rest of our analysis, we denote by $'$ the $\tau$-derivative. In the new time, $\varphi'(\tau) = \omega$ and so $\varphi(\tau) = \omega \tau$. Hence it suffices to solve
    \begin{equation}    \label{reducedsystem}           
        \begin{aligned}
        \tilde{r}'&=\frac{\varepsilon\beta_2\tilde{r}\tilde{z} + \varepsilon^2\left(\tilde{\mu}\gamma_3(\omega\tau)\tilde{r} + \tilde{\mu}\gamma_4(\omega\tau)\tilde{z} + \beta_3\tilde{r}^3 + \beta_4\tilde{r}\tilde{z}^2\right) 
 + O(\varepsilon^3)}{1 + \varepsilon\beta_1 \tilde{z} + O(\varepsilon^2)},\\
        \tilde{z}'&=\frac{\varepsilon\left(\tilde{\mu}\gamma_5 + \beta_5 \tilde{r}^2\right)+\varepsilon^2\left(\tilde{\mu}\gamma_6(\omega\tau)\tilde{r} + \tilde{\mu}\gamma_7\tilde{z} + \beta_6\tilde{r}^2\tilde{z}\right) + O(\varepsilon^3)}{1 + \varepsilon\beta_1 \tilde{z} + O(\varepsilon^2)}.
        \end{aligned}
    \end{equation}
Expanding \eqref{reducedsystem} in $\varepsilon > 0$ leads to 
\begin{equation}\label{redeq2}
     \begin{aligned}
        \tilde{r}'&=\varepsilon\beta_2\tilde{r}\tilde{z} + \varepsilon^2\left(\tilde{\mu}\gamma_3 \tilde{r}+\tilde{\mu}\gamma_4 \tilde{z} + \beta_3\tilde{r}^3 - (\beta_1\beta_2 - \beta_4) \tilde{r}\tilde{z}^2\right) +O(\varepsilon^3),\\
         \tilde{z}'&=\varepsilon\left(\tilde{\mu}\gamma_5 + \beta_5 \tilde{r}^2\right) + \varepsilon^2\left(\tilde{\mu}\gamma_6\tilde{r} - \tilde{\mu}(\beta_1\gamma_5 -  \gamma_7) \tilde{z} - (\beta_1\beta_5 - \beta_6) \tilde{r}^2\tilde{z}\right) + O(\varepsilon^3),
    \end{aligned}
\end{equation}
where we have omitted the $\tau$-dependence of $\gamma_3,\gamma_4,$ and $\gamma_6$ for simplicity of notation.

\textbf{Step 3 (existence of periodic orbits):} The first-order truncated ODEs of \eqref{redeq2} in $\varepsilon$ read
\begin{equation}\label{Eq:1jet}
     \begin{aligned}
         \tilde{r}'&=\varepsilon\beta_2\tilde{r}\tilde{z},\\
         \tilde{z}'&=\varepsilon\left(\tilde{\mu}\gamma_5 + \beta_5 \tilde{r}^2\right),
    \end{aligned}
\end{equation}
and have a nonzero equilibrium
\begin{equation}\label{equilibrium}
    (\tilde{r}_0,\tilde{z}_0):=\left(\sqrt{\frac{-\tilde{\mu}\gamma_5}{\beta_5}},0\right) \quad \mbox{if and only if} \quad \tilde{\mu}\beta_5\gamma_5<0.
\end{equation}
The Jacobian matrix of \eqref{Eq:1jet} at $(\tilde{r}_0, \tilde{z}_0)$ is

\begin{equation}\label{first-order}
    \varepsilon A_1:=\varepsilon\left(\begin{array}{cc}
               0  &  \beta_2 \tilde{r}_0 \\
               2\beta_5 \tilde{r}_0  &  0
            \end{array}\right).
\end{equation}
Fixing $0 < \tilde{r}_0 < 1$, since $A_1$ is invertible, averaging theory \cite[Theorem 6.3.2]{SanVer} ensures a continuation of $(\tilde{r}_0, \tilde{z}_0)$  for sufficiently small $0 < \varepsilon \ll 1$. More precisely, there exists a smooth $\varepsilon$-family of periodic solutions of the full ODEs \eqref{redeq2}, denoted by
    \begin{equation}\label{persol}
        (\tilde{r}(\tau;\varepsilon), \tilde{z}(\tau;\varepsilon)):=(\tilde{r}_0,\tilde{z}_0) + \varepsilon (\tilde{R}(\tau;\varepsilon),\tilde{Z}(\tau;\varepsilon)),
    \end{equation}
and moreover, the minimal period is $2\pi/\omega$. 

To see the $\varepsilon$-dependence of all variables explicitly, we now recover the original variables (i.e., without the tilde) via the rescaling \eqref{kappaScaling}. Notice that the choice $0 < \tilde{r}_0 < 1$ ensures that the equilibrium $(\tilde{r}_0, \tilde{z}_0)$  in \eqref{equilibrium} resides in the region of validity $\Phi(U_1)$ of the expansion \eqref{redeq2}, yielding an upper bound $\tilde{\mu}_0$ for $\lvert \tilde{\mu} \rvert$ given by
\begin{equation}
0 < \lvert \tilde{\mu} \rvert <\tilde{\mu}_0<\left \lvert \frac{\beta_5}{\gamma_5}\right \rvert.
\end{equation}
The quantity $\mu_0$ in Theorem \ref{T:Bifurcation} reads
\begin{equation}
    \mu_0 = \varepsilon^2 \left \lvert \frac{\beta_5}{\gamma_5}\right \rvert.
\end{equation}
The direction of bifurcation of periodic solutions is determined by $\mathrm{sign}(\tilde{\mu}) = \mathrm{sign}(\mu)$, which is equal to $-\mathrm{sign}(\beta_5\gamma_5)$ by \eqref{equilibrium}, and thus  \eqref{Eq:Direction} is proved.

As $\varepsilon \searrow 0$, the original variables have the asymptotics

\begin{equation}
    r(t;\varepsilon) = \varepsilon \tilde{r}_0 + O(\varepsilon^2) = \sqrt{\frac{-\mu\beta_5}{\gamma_5}} + O(\varepsilon^2), \quad z(t; \varepsilon) = O(\varepsilon^2), \quad \text{and} \quad \mu = O(\varepsilon^2),
\end{equation}
where the $O(\varepsilon^2)$-terms depend on $t$ periodically. Therefore, $r(t;\varepsilon) = O(\sqrt{\lvert \mu \rvert})$ and $z(t; \varepsilon) = O(\mu)$ by \eqref{equilibrium} and \eqref{persol}. The limit \eqref{Eq:LimitPeriod} of minimal period follows from $\lim_{\varepsilon \searrow 0}\dot\varphi(t;\varepsilon) =  \omega$, since \eqref{persol} is also a periodic solution of the ODEs \eqref{perturbedform2}.

\textbf{Step 4 (stability criteria for periodic orbits):} We distinguish between the hyperbolic case and the elliptic case.

In the hyperbolic case, i.e., $\xi= \mathrm{sign}(\beta_2\beta_5) = 1$, the matrix $A_1$ in \eqref{first-order} is hyperbolic because it has a positive eigenvalue and a negative eigenvalue. Thus, for sufficiently small $0 < \varepsilon \ll 1$, every periodic solution of the original ODEs \eqref{StandardSystem} is exponentially unstable; see \cite[Theorem 6.3.3]{SanVer}. Moreover, as we recover the angle variable $\varphi$, the periodic solutions in the original ODEs \eqref{StandardSystem} are associated with a two-dimensional unstable manifold.

In the elliptic case, i.e., $\xi= \mathrm{sign}(\beta_2\beta_5) = -1$, the matrix $A_1$ defined in \eqref{first-order} has two purely imaginary simple eigenvalues and thus it is not hyperbolic. Hence we cannot conclude the local stability properties solely from the first-order truncated ODEs \eqref{Eq:1jet}; see \cite[Chapter 6]{SanVer}. We thereby analyze the second-order truncated ODEs of \eqref{redeq2}.

We first approximate the Floquet multipliers associated with the periodic solution \eqref{persol}. To this end, differentiating \eqref{redeq2} along \eqref{persol} yields the linear equations

    \begin{equation}\label{lineq}
    \begin{aligned}
        u'&=\varepsilon\beta_2\tilde{r}_0v + \varepsilon^2 \left(\left(3\beta_3\tilde{r}_0^2 + \tilde{\mu}\gamma_3 + \beta_2 \tilde{Z}\right)u + (\tilde{\mu}\gamma_4 + \beta_2 \tilde{R})v\right) + O(\varepsilon^3),\\
        v'&= 2\varepsilon \beta_5\tilde{r}_0u + \varepsilon^2\left((\tilde{\mu}\gamma_6 + 2\beta_5 \tilde{R})u - \left(\beta_1\beta_5 \tilde{r}_0^2  - \beta_6\tilde{r}_0^2 +\tilde{\mu}\beta_1\gamma_5- \tilde{\mu}\gamma_7\right)v\right) + O(\varepsilon^3).
        \end{aligned}
    \end{equation}
Since the periodic solution \eqref{persol} is $2\pi/\omega$-periodic, by the variation-of-constants formula, the Floquet multipliers are $O(\varepsilon^2)$-close to the eigenvalues of $\mathrm{e}^{2\pi\varepsilon A_1/\omega}$, where $A_1$ is defined in \eqref{first-order}.

Choosing sufficiently small $0 < \varepsilon \ll 1$ ensures that the Floquet multipliers of the periodic solution \eqref{persol} are $\{\kappa,\overline{\kappa}\}$ for some $\kappa\in \mathbb{C}\setminus\mathbb{R}$. The matrix that collects all $\varepsilon^2$-terms in the Jacobian matrix of \eqref{lineq} at the periodic solution \eqref{persol} reads 

    \begin{equation}\label{second-order}
           A_2(\tau ):=\left(\begin{array}{cc}
           3\beta_3\tilde{r}_0^2+\tilde{\mu}\gamma_3(\omega\tau) + \beta_2 \tilde{Z}(\tau;\varepsilon ) & \tilde{\mu}\gamma_4(\omega\tau ) + \beta_2 \tilde{R}(\tau;\varepsilon)\\
           \tilde{\mu} \gamma_6(\omega \tau) + 2 \beta_5 \tilde{R}(\tau;\varepsilon )  & -\beta_1\beta_5 \tilde{r}_0^2+ \beta_6\tilde{r}_0^2 -\tilde{\mu}\beta_1\gamma_5+\tilde{\mu}\gamma_7
        \end{array}
        \right).
    \end{equation}
By Liouville's formula, we obtain the following relation of Floquet multipliers:
    \begin{equation}
        \kappa \overline{\kappa} = \lvert \kappa \rvert^2 =\mathrm{exp}\left(\varepsilon^2\int_0^{\frac{2\pi}{\omega}}\mathrm{tr}(A_2(\tau))\,\mathrm{d}\tau+O(\varepsilon^3)\right),
    \end{equation}
where $\mathrm{tr}(\cdot)$ denotes the trace of matrices. Hence the sign of $\int_0^{2\pi/\omega}\mathrm{tr}(A_2(\tau ))\,\mathrm{d}\tau$ determines the local stability properties of the periodic solution. Indeed,

    \begin{equation} \label{nontrivial-integral}
        \begin{aligned}
            \int_0^{\frac{2\pi}{\omega}}\mathrm{tr}(A_2(\tau ))\,\mathrm{d}\tau=\int_0^{\frac{2\pi}{\omega}}\big(&  -\beta_1\beta_5 \tilde{r}_0^2 + 3\beta_3\tilde{r}_0^2 + \beta_6\tilde{r}_0^2 
           \\
            &- \tilde{\mu}\beta_1\gamma_5 +\tilde{\mu}\gamma_3(\omega\tau)+\tilde{\mu}\gamma_7 + \beta_2\tilde{Z}(\tau;\varepsilon)\big)\,\mathrm{d}\tau,
    \end{aligned}
    \end{equation}
in which the only nonconstant integrand is $\tilde{\mu}\gamma_3(\omega\tau)+\beta_2 \tilde{Z}(\tau;\varepsilon)$.

To compute \eqref{nontrivial-integral}, we substitute the periodic solution \eqref{persol} into the full ODEs \eqref{redeq2} and obtain
    \begin{equation} \label{R-equation}
        \tilde{R}'(\tau;\varepsilon)= \varepsilon \left(\beta_3\tilde{r}_0^3 + \tilde{\mu}\gamma_3(\omega\tau)\tilde{r}_0 + \beta_2 \tilde{r}_0 \tilde{Z}(\tau;\varepsilon) \right) + O(\varepsilon^2).
    \end{equation}
Since $\tilde{R}(\tau;\varepsilon)$ is $2\pi/\omega$-periodic in $\tau$, integrating \eqref{R-equation} over $[0, 2\pi/\omega]$ yields
    \begin{equation} \label{nontrivial-calculation}
        \int_0^{\frac{2\pi}{\omega}}
        \tilde{\mu} \gamma_3(\omega\tau)+ \beta_2 \tilde{Z}(\tau;\varepsilon) \,\mathrm{d}\tau=\frac{-2\pi\beta_3 \tilde{r}_0^2}{\omega}+O(\varepsilon).
    \end{equation}
Recalling $\omega>0$, as $0 < \varepsilon \ll 1$ we use \eqref{nontrivial-calculation} and substitute $\tilde{r}_0 = \sqrt{-\tilde{\mu} \gamma_5/\beta_5}$ to derive
    \begin{equation}\label{signeq}
        \begin{aligned}
        \mathrm{sign}(\log \lvert \kappa \rvert)&=\mathrm{sign}\left(\int_0^{\frac{2\pi}{\omega}}\mathrm{tr}(A_2(\tau))\, \mathrm{d}\tau\right)\\
        &=\mathrm{sign}\left(-\beta_1\beta_5\tilde{r}_0^2+2\beta_3\tilde{r}_0^2+ \beta_6\tilde{r}_0^2 -\tilde{\mu}\beta_1\gamma_5+\tilde{\mu}\gamma_7\right)\\
        &=\mathrm{sign}\left(-\frac{2\tilde{\mu}\beta_3\gamma_5}{\beta_5}-\frac{\tilde{\mu}\beta_6\gamma_5}{\beta_5}+\tilde{\mu}\gamma_7\right).
        \end{aligned}
    \end{equation}
By \eqref{equilibrium}, the periodic solutions bifurcate for $\mathrm{sign}(\tilde{\mu}\beta_5\gamma_5)=-1$. Thus, we multiply the right-hand-side of \eqref{signeq} by $-\beta_5\gamma_5/\tilde{\mu}$, note $\mathrm{sign}(-\beta_5 \gamma_5/\tilde{\mu}) = 1$, and conclude 
    \begin{equation}
        \mathrm{sign}(\log \lvert \kappa \rvert)=\mathrm{sign}\left(2\beta_3\gamma_5^2 -\beta_5\gamma_5\gamma_7 + \beta_6\gamma_5^2\right).
    \end{equation}
Hence $\lvert \kappa \rvert < 1$ (resp., $\lvert \kappa \rvert > 1$) if and only if $2\beta_3\gamma_5^2-\beta_5\gamma_5\gamma_7 +\beta_6\gamma_5^2 < 0$ (resp., $> 0$), which proves the stability criteria \eqref{Eq:Stability}--\eqref{Eq:Instability}. 
\end{proof}

\begin{remark}(Motivation for Liouville's Lemma)
As noted in Step 4 above, the stability of periodic orbits near an elliptic Hopf point cannot be determined from the first-order truncated ODEs \eqref{first-order}. It is tempting to perform second-order averaging and consider the $\tau$-average of $\varepsilon A_1 + \varepsilon^2 A_2(\tau)$, where $A_2(\tau)$ is given by \eqref{second-order}. However, this analysis is not sufficient to conclude stability, because the characteristic equation for eigenvalues associated to the second-order averaged system depends on $\varepsilon$ in a nonlinear manner. In particular, the $\varepsilon$-expansion of the eigenvalues may contain terms of order higher than two. Hence higher-order averaging is required. Instead, our approach has the advantage that it only uses the trace, and thus it prevents the need for a higher-order expansion.
\end{remark}

\begin{proof}[\textbf{Proof of Lemma \ref{lemma-all-elliptic}}.]
We transform the predator-prey system \eqref{unfolding} with $\mu = 0$ into the ODEs \eqref{StandardSystem} satisfying the assumptions \textbf{A1}--\textbf{A3}, by applying an affine linear change of coordinates such that Hopf point $X_H$ is the origin and the Jacobian matrix associated with $X_H$ is in the Jordan form. 

In the admissible parameter region \eqref{ParRegion} we have 
\begin{equation} \label{positive-parameters}
\ell_1:=1-2\lambda-\alpha_1> 0 \quad \mbox{and} \quad \ell_2 :=2 \lambda + \alpha_2-1 > 0.
\end{equation}
Then the frequency satisfies 
     \begin{equation}
    \omega=\sqrt{\frac{\lambda(\delta_1 \ell_2+\delta_2 \ell_1)}{\ell_1+\ell_2}}.
     \end{equation}
Recall that $\xi = \mathrm{sign}(\beta_2 \beta_5)$ is
the discriminant for Hopf bifurcations without parameter; see \eqref{discriminant}. From the formulas in Appendix \ref{Appendix} we obtain
\begin{equation}
         \beta_2 =\frac{-\lambda(\ell_1+\ell_2)}{2(\lambda + \alpha_1)(\lambda + \alpha_2)^2} < 0 \quad \mbox{and} \quad \beta_5 =\frac{\lambda\delta_1\delta_2\ell_1\ell_2}{2(\lambda+\alpha_1)\omega^2(\ell_1+\ell_2)} > 0.
     \end{equation}
Hence the assumption \textbf{A4} is fulfilled with $\xi = -1$, i.e., the Hopf point $X_H$ is elliptic.
\end{proof}

\begin{proof}[\textbf{Proof of Theorem \ref{T:StableCoexistence}}.]
As in the proof of Lemma \ref{lemma-all-elliptic}, we apply an affine linear change of coordinates that transforms the predator-prey system \eqref{unfolding} into the ODEs \eqref{StandardSystem} satisfying the assumptions \textbf{A1}--\textbf{A4}. Then we modify the coordinates via \eqref{Eq:ParamShift} to fulfill \textbf{A5}.

We use computer algebra software (e.g., Maple or Maxima) to obtain the coefficients of the perturbed cylindrical form \eqref{perturbedform}; see Lemma \ref{L:PertNormalform}. We set $\ell_1:=1-2\lambda-\alpha_1 > 0$ and $\ell_2:=2 \lambda + \alpha_2-1 > 0$. Then the coefficients for the first-order truncated ODEs \eqref{Eq:1jet} are 
     \begin{equation} \label{PP-betas}
     \begin{aligned}
         \beta_2&=\frac{-\lambda(\ell_1+\ell_2)}{2(\lambda + \alpha_1)(\lambda + \alpha_2)^2} < 0,\\
         \beta_5&=\frac{\lambda\delta_1\delta_2\ell_1\ell_2}{2(\lambda+\alpha_1)\omega^2(\ell_1+\ell_2)} > 0,\\
         \gamma_5&=\frac{-\lambda(\lambda+\alpha_2)\delta_1\delta_2\ell_1\ell_2}{\omega^2(\ell_1+\ell_2)^2} < 0;
     \end{aligned}
     \end{equation}
see the formulas in Appendix \ref{Appendix}. Since $\mathrm{sign}(\mu) = -\mathrm{sign}(\beta_5 \gamma_5) = 1$ by \eqref{Eq:Direction}, the branch of periodic orbits obtained in Theorem \ref{T:Bifurcation} emanates for $\mu > 0$. 

To verify the stability inequality \eqref{Eq:Stability}, we compute the coefficients in the perturbed cylindrical form \eqref{perturbedform} from Appendix \ref{Appendix}:
     \begin{equation} \label{PP-gammas}
     \begin{aligned}
         \beta_3&=\frac{\lambda ((\lambda + \alpha_1)\delta_1\ell_2H_1 - (\lambda + \alpha_2)\delta_2\ell_1H_2)}{8(\lambda+\alpha_1)^2(\lambda+\alpha_2^2)\omega^2(\ell_1+\ell_2)}+\frac{\lambda^2\delta_1\delta_2\ell_1\ell_2((\lambda+\alpha_1)\delta_2-(\lambda+\alpha_2)\delta_1)}{4(\lambda+\alpha_1)^2(\lambda+\alpha_2)^2\omega^4(\ell_1+\ell_2)},
         \\         \beta_6&=\frac{\lambda^2(\lambda+\alpha_1+\ell_1)\delta_1\delta_2(\delta_1\ell_2-\delta_2\ell_1)}{2(\lambda+\alpha_1)^2(\lambda+\alpha_2)^2\omega^4},\\
         \gamma_7&=\frac{-\lambda^2\delta_1\delta_2\left((\lambda+\alpha_1)\delta_1\ell_2^2-(\lambda+\alpha_2)\delta_2\ell_1^2\right)}{(\lambda+\alpha_1)(\lambda+\alpha_2)\omega^4(\ell_1+\ell_2)^2}.
     \end{aligned}
     \end{equation}
Here $H_1$ and $H_2$ are the following quadratic polynomials:
     \begin{equation}
     \begin{aligned}
         H_1(\lambda,\alpha_1,\alpha_2)&=- \lambda + 2\alpha_2 + \lambda\alpha_1 - 8\lambda\alpha_2 - 2\alpha_1\alpha_2  - 2\alpha_2^2,\\
         H_2(\lambda,\alpha_1,\alpha_2)&=\lambda - 2\alpha_1
         + 8\lambda\alpha_1 - \lambda\alpha_2 + 2\alpha_1\alpha_2 + 2\alpha_1^2.
     \end{aligned}
     \end{equation}
Notice that both $H_1$ and $H_2$ are independent of $\delta_1$ and $\delta_2$, motivating the viewpoint of tetrahedra on the admissible parameter region \eqref{ParRegion}; see Remark \ref{remark-tetrahedra}. Moreover, note $H_2(\lambda,\alpha_1,\alpha_2)=-H_1(\lambda,\alpha_2,\alpha_1)$. 
By using the coefficients in \eqref{PP-betas}--\eqref{PP-gammas}, we verify that the stability inequality \eqref{Eq:Stability} is equivalent to 
     \begin{equation}\label{Ineq:Stab}
     (\lambda+\alpha_1) \delta_1 \ell_2H_1-(\lambda+\alpha_2)\delta_2\ell_1 H_2<0.
     \end{equation}

We perform the subsequent analysis in the admissible parameter region \eqref{ParRegion}. Then $\ell_1>0$ and $\ell_2>0$, and so it suffices to prove $H_1(\lambda, \alpha_1, \alpha_2) <0$ and $H_2(\lambda, \alpha_1, \alpha_2)>0$. 

The key feature is that $H_1$ is affine linear in $\alpha_1$ for $0 < \alpha_1 < 1 - 2 \lambda$, and thus $H_1(\lambda, \alpha_1, \alpha_2)$ has a fixed sign if
both $H_1(\lambda, 0, \alpha_2)$ and $H_1(\lambda, 1 - 2\lambda, \alpha_2)$ share the same sign. Indeed
     \begin{equation}
         H_1(\lambda, 0,\alpha_2)=  - \lambda+(2-8\lambda)\alpha_2 -2\alpha_2^2 \quad \text{and}\quad H_1(\lambda, 1-2\lambda,\alpha_2)=-2(\lambda + \alpha_2)^2.
     \end{equation}
Clearly, $H_1(\lambda, 1-2\lambda,\alpha_2)<0$. We observe that the graph of $H_1(\lambda, 0, \alpha_2)$ in $\alpha_2$ is a parabola facing downwards. Since $\alpha_2 > 1 - 2 \lambda > 0$ and so $H_1(\lambda, 0, 1 - 2 \lambda) = 8 \lambda^2 - 5 \lambda < 0$, we know $H_1(\lambda, 0, \alpha_2) < 0$. Hence $H_1(\lambda, \alpha_1, \alpha_2) < 0$ in the region \eqref{ParRegion}. 

Analogously, $H_2$ is affine linear in $\alpha_2$ for $1 -2 \lambda  < \alpha_2 <1$. Since
     \begin{equation}
         H_2(\lambda, \alpha_1,1-2\lambda)=2(\lambda + \alpha_1)^2 > 0 \quad \text{and}\quad H_2(\lambda, \alpha_1,1)=8\lambda\alpha_1 + 2\alpha_1^2 > 0,
     \end{equation}
we know $H_2(\lambda, \alpha_1, \alpha_2) > 0$ in the region \eqref{ParRegion}. Consequently, \eqref{Ineq:Stab} holds in the admissible parameter region \eqref{ParRegion}. The proof is complete. 
 \end{proof}

It is challenging to verify the stability inequality \eqref{Eq:Stability} for all parameter values in the admissible region \eqref{ParRegion}, even with computer algebra software. There are two reasons. First, the stability coefficient $2\beta_3\gamma_5^2-\beta_5\gamma_5\gamma_7+\beta_6\gamma_5^2$ in \eqref{Eq:Stability} is a rational function of the parameters and contains polynomials of degree ten in five variables. Second, the region \eqref{ParRegion} involving $\delta_1, \delta_2 > 0$ is unbounded.

Our first approach to explore the stability inequality \eqref{Eq:Stability} was a numerical search of maximizers within the region \eqref{ParRegion} by Newton's method, for bounded values of $\delta_1$ and $\delta_2$. Taking random initial conditions, we saw that the stability coefficient remained negative in the region \eqref{ParRegion} and the maximizers accumulated to the boundary of \eqref{ParRegion}.

Motivated by our numerical exploration, we proceeded to prove Theorem \ref{T:StableCoexistence}. The key observation is that the stability inequality \eqref{Eq:Stability} simplifies, as we identify the new parameters $\ell_1 := 1 - 2\lambda - \alpha_1$ and $\ell_2 := 2\lambda + \alpha_2 - 1$. Note that $\ell_j$ for $j = 1, 2$ is meaningful in our ecological system \eqref{unfolding}, because $\ell_j = 0$ is the value of a classical Hopf bifurcation on the boundary plane $Q_j$; see \eqref{boundary-Hopf-value}. Then by carefully examining \eqref{Eq:Stability}, we derived the simplified inequality \eqref{Ineq:Stab}, making the remaining part of the proof manageable.

\appendix
\section{Formulas of Coefficients in the cylindrical forms}\label{Appendix}

We list the formulas of the coefficients $\beta_j$ and $\gamma_j$ in the cylindrical forms in Lemma \ref{L:UnpertNormalForm} and Lemma \ref{L:PertNormalform}. We design notations in the way that the coefficients $\beta_j$ correspond to the unperturbed cylindrical form (i.e., $\mu = 0$), while the other coefficients $\gamma_j$ correspond to the perturbed cylindrical form (i.e., $\mu \neq 0$) and so they involve the $\mu$-derivatives. According to Theorem \ref{T:Bifurcation}, we highlight the role of the coefficients as follows.
\begin{itemize}
\item The direction of bifurcation is determined by both $\beta_5$ and $\gamma_5$.
\item The types of bifurcation, \textbf{Type-H} and \textbf{Type-E}, are determined by both $\beta_2$ and $\beta_5$.
\item \textbf{Type-ES} and \textbf{Type-EU} are determined by all $\beta_3$, $\beta_5$, $\beta_6$, $\gamma_5$, and $\gamma_7$.
\end{itemize}
\textit{All derivatives below are evaluated at $(y,z;\mu) = (0,0;0)$.}
\begin{equation}\nonumber
        \begin{aligned}
            \beta_1&= \frac{1}{2\omega}\left(\partial_{y_1}\partial_{z}f^{y_2}-\partial_{y_2}\partial_{z}f^{y_1}\right),\\                  \beta_2&= \frac{1}{2}\partial_z\mathrm{div}_yf^y,\\
             \beta_3&= \frac{1}{16}\left(\partial_{y_1}\Delta_yf^{y_1}+\partial_{y_2}\Delta_yf^{y_2}\right)+\frac{1}{16\omega}\left(\partial_{y_1}\partial_{y_2}f^{y_1}\Delta_yf^{y_1}-\partial_{y_1}\partial_{y_2}f^{y_2}\Delta_yf^{y_2}\right)\\   &\quad+\frac{1}{16\omega}\left(\partial_{y_2}^2f^{y_1}\partial_{y_2}^2f^{y_2}-\partial_{y_1}^2f^{y_1}\partial_{y_1}^2f^{y_2}\right) \\ &
            \quad +\frac{1}{16\omega}\left(\partial_{y_2}\partial_zf^{y_2}-\partial_{y_1}\partial_zf^{y_1}\right)\partial_{y_1}\partial_{y_2}f^z\\            & \quad+\frac{1}{32\omega}\left(\partial_{y_2}\partial_zf^{y_1}+\partial_{y_1}\partial_zf^{y_2}\right)\left(\partial_{y_1}^2f^z-\partial_{y_2}^2f^z\right), \\
             \beta_4&= \frac{1}{4}\partial_z^2\mathrm{div}_yf^y,\\
            \beta_5&= \frac{1}{4}\Delta_y f^z,\\
            \beta_6&=\frac{1}{4}\partial_z\Delta_yf^z+\frac{1}{4\omega}\left(\partial_{y_2}\partial_zf^z\Delta_yf^{y_1}-\partial_{y_1}\partial_zf^z\Delta_yf^{y_2}\right)\\
            &\quad +\frac{1}{4\omega}\left(\partial_{y_1}\partial_zf^{y_1}-\partial_{y_2}\partial_zf^{y_2}\right)\partial_{y_1}\partial_{y_2}f^z\\
            &\quad + \frac{1}{8\omega}\left(\partial_{y_2}\partial_zf^{y_1}+\partial_{y_1}\partial_zf^{y_2}\right)\left(\partial_{y_2}^2f^z-\partial_{y_1}^2f^z\right),
            \\
        \end{aligned}
\end{equation}

\begin{equation}\nonumber
        \begin{aligned}
              \gamma_1(\varphi)&=-\partial_\mu\partial_{y_2}f^{y_1}\sin^2(\varphi) +(\partial_\mu\partial_{y_2}f^{y_2}-\partial_\mu\partial_{y_1}f^{y_1})\sin(\varphi)\cos(\varphi)\\
            & \quad + \partial_\mu\partial_{y_1}f^{y_2}\cos^2(\varphi) - \frac{1}{2\omega}\left(\partial_{y_1}\partial_zf^{y_1} - \partial_{y_2}\partial_zf^{y_2}\right)\partial_\mu f^z\cos^2(\varphi)\\
            &\quad - \frac{1}{2\omega} \left(\partial_{y_2}\partial_zf^{y_1} + \partial_{y_1}\partial_zf^{y_2}\right)\partial_\mu f^z\sin(\varphi)\cos(\varphi)\\
            &\quad -\frac{1}{2\omega}\left( \pi \partial_{y_2}\partial_z f^{y_1} - \pi \partial_{y_1}\partial_z f^{y_2} - \frac{1}{2}\partial_{y_1}\partial_z f^{y_1} + \frac{1}{2}\partial_{y_2}\partial_z f^{y_2}\right)\partial_\mu f^z,
            \\
            \gamma_2(\varphi)&= -\partial_\mu\partial_{z}f^{y_1}\sin(\varphi)+\partial_\mu\partial_{z}f^{y_2}\cos(\varphi),\\
            \gamma_3(\varphi)&=\partial_\mu\partial_{y_2}f^{y_2}\sin^2(\varphi)+(\partial_\mu\partial_{y_2}f^{y_1}+\partial_\mu\partial_{y_1}f^{y_2})\sin(\varphi)\cos(\varphi) \\ & \quad + \partial_\mu\partial_{y_1}f^{y_1}\cos^2(\varphi) + \frac{1}{2\omega}\left(\partial_{y_2}\partial_zf^{y_1} + \partial_{y_1}\partial_zf^{y_2}\right)\partial_\mu f^z\cos^2(\varphi)\\
            &\quad - \frac{1}{2\omega} \left(\partial_{y_1}\partial_zf^{y_1} - \partial_{y_2}\partial_zf^{y_2}\right)\partial_\mu f^z\sin(\varphi)\cos(\varphi)
            \\
            &\quad + \frac{1}{2\omega}\left( \pi \partial_{y_1}\partial_z f^{y_1} + \pi \partial_{y_2}\partial_z f^{y_2} - \frac{1}{2}\partial_{y_2}\partial_z f^{y_1} - \frac{1}{2}\partial_{y_1}\partial_z f^{y_2}\right)\partial_\mu f^z,\\
            \gamma_4(\varphi)&= \partial_\mu\partial_{z}f^{y_2}\sin(\varphi)+ \partial_\mu\partial_{z}f^{y_1}\cos(\varphi),\\
            \gamma_5&= \partial_\mu f^z,\\
            \gamma_6(\varphi)&= \left(\partial_\mu\partial_{y_2}f^{z}-\frac{1}{\omega}\partial_{y_1}\partial_z f^z\partial_\mu f^z\right)\sin(\varphi)\\
            & \quad + \left(\partial_\mu\partial_{y_1}f^{z} + \frac{1}{\omega}\partial_{y_2}\partial_z f^z\partial_\mu f^z\right)\cos(\varphi),\\
            \gamma_7&=\partial_\mu\partial_{z}f^{z}.\\
            \end{aligned}
        \end{equation}

\textbf{Acknowledgement.}
We are grateful to Sze-Bi Hsu for the suggestion of studying stable periodic coexistence and to Bernold Fiedler for many inspiring discussions. 

\textbf{Funding.} A. L.-N. has been supported by NSTC grant 113-2123-M-002-009. P. L. has been supported by Marie Skłodowska--Curie Actions, UNA4CAREER H2020 Cofund, 847635, with the project DYNCOSMOS. N.V. has been supported by the DFG (German Research Society), project n. 512355535. J.-Y. D. has been supported by NSTC grant 113-2628-M-007-005-MY4.




\end{document}